\documentclass{article}

\title{$R\text{-}\mathrm{Mod}$-enriched categories are left $\un{R}$-module objects of $Cat(\Ab)$ and $Cat(\Ab)$-enriched functors }
\author{Matteo Doni\footnote{Università degli Studi di Milano, Milan. Email address: matteo.doni@unimi.it}}
\input{CategoricalArticle.sty}
\begin{document}
\newlength{\myindent} 
\setlength{\myindent}{\parindent}
\parindent 0em 
\maketitle
\begin{abstract}
We establish the feasibility of investigating the theory of $R\text{-}\mathrm{Mod}$-enriched categories, for any unitary ring  $R$, through the framework of $\Ab$-enriched category theory. In particular, we prove that the category of $R$-$\mathrm{Mod}$-enriched categories, $Cat(R$-$\mathrm{Mod})$, the category of $\un{R}$-modules inside $Cat(\Ab)$, $\mathrm{LMod}_{\un{R}}(Cat(\Ab))$, and the category of $Cat(\Ab)$-enriched functors, $Fun^{Cat(\Ab)}(\2un{R},Cat(\Ab))$, are equivalent.
\end{abstract}
\setcounter{tocdepth}{1}
\tableofcontents

\section{Introduction}
The theory of enriched categories is a natural extension of the theory of categories. Some categories have hom-sets that possess structures beyond simple sets. For instance, the category of (nice) topological spaces $\mathrm{Top}$ has hom-sets that are topological spaces with the compact-open topology. The category of $\R$-vector spaces has hom-sets that are $\R$-vector.
The category of small categories $\mathrm{Cat}$ has hom-sets that are small categories: for two categories $X$ and $Y$, the hom-set is the category whose objects are functors and whose morphisms are natural transformations.

Enriched category theory is such a natural extension of category theory that it has always been used to study categorical settings; it is unthinkable to have a theory of categories without natural transformations. Indeed, the definition of natural transformation is in the first reference book of category theory \cite{MacLane}.  

Enriched category theory is ubiquitous in mathematics.
In general, hom-sets can be objects in any monoidal category $(\V, -\otimes_{\V}-,\I_{\V})$, and these are called \textit{hom-$\V$-objects} or, simply, \textit{hom-objects}.
An enriched category with hom-$\V$-objects is called a $\V$-enriched category; see \Cref{defEnichedCategory}.
Moreover, when comparing two $\V$-enriched categories, one must account for the additional structure of the hom-objects, which naturally leads to the definition of $\V$-enriched functors; see \Cref{defEnrichedFunctor}.   

By the notation $\mathrm{Cat}(\V)$, we denote the category of small $\V$-categories whose morphisms are $\V$-enriched functors; see \Cref{notCatV}. 

The motivations presented so far may not suffice to justify the effort of developing this theory. The necessity for this theory arises from its ability to encompass certain homotopic and (co)homological phenomena that ordinary category theory cannot. For example, homotopy (co)limits can be computed as enriched (co)limits \cite[\S 6.6]{riehl2014categorical}. 
It is not coincidental that $\infty$-categories are $\mathrm{Top}$-enriched categories (with certain properties) or that, in derived algebraic geometry, where homology and cohomology relations are prevalent, the theory of $Ch(\Z)$-enriched categories, where $Ch(\Z)$ is the category of chain complexes (commonly known as dg-category theory), is widely used.

Considering hom-objects instead of hom-sets can be daunting and might seem to complicate matters. However, this complication can sometimes simplify problem-solving, which is why it is sometimes intentionally induced. This is precisely the idea that leads to the definition of stable homotopy theory \cite[\S 1.4]{HA}.

In this paper, we are not concerned with homotopy, (co)homology, or other related topics.
We focus on classical linear categories, specifically $\Ab$-enriched categories, where $\Ab$ is the category of abelian groups, or categories whose hom-objects lie within certain extensions of $\Ab$, i.e., objects in the category of $R$-modules $R$-$\mathrm{Mod}$, where $R$ is a unitary ring.    

The goal of this paper is to show that it is possible to study the category of $R$-$\mathrm{Mod}$-enriched categories, $\mathrm{Cat}(R$-$\mathrm{Mod})$, within the category of $\Ab$-enriched categories, $\mathrm{Cat}(\Ab)$.

This result is expected because it is classical to use the extension-restriction adjunction:
\begin{equation}
    \label{eqextentionrestri}
   -\otimes_{\Z} R :\Ab\leftrightarrows R-\mathrm{Mod}:U; 
\end{equation}
to retrieve information about $R$-$\mathrm{Mod}$ from $\Ab$; see \cite[Example 5.5.7]{RiehlCTIC}.

Unlike $R$-$\mathrm{Mod}$-enriched categories, $\Ab$-enrichments have interesting properties because the $\Ab$-enrichment can be (at least partially) internalized. For example, in abelian or additive category theory, the operations in the hom-objects of an $\Ab$-enriched category are linked to its (co)limits.

Using the functor $U$ in \eqref{eqextentionrestri}, it is possible to disregard the action of $R$ and consider the hom-objects of an $R$-$\mathrm{Mod}$-enriched category simply as hom-$\Ab$-objects (this is often referred to as base change or change of enrichment).   

Therefore, to achieve our goal, we need to find a way to recover the forgotten action of $R$ within the category $\mathrm{Cat}(\Ab)$. In this paper, we present two methods to accomplish this.

The first idea is very natural.
We need to internalize an action of a ring $R$ within $\mathrm{Cat}(\Ab)$. The first step is to internalize the ring $R$. 
Any ring $R$ can be considered as an $\Ab$-enriched category $\underline{R}$ with only one object $\{*\}$, whose only non-trivial hom-object is $\Hom(*,*)=R$, and whose composition is the multiplication in $R$. This internalization is widely used; for example, dg-algebras are monoids (where monoids are equivalent to rings) internalized in a $Ch(\Z)$-enriched category.
The second step is to define the action of $\underline{R}$ within $\mathrm{Cat}(\Ab)$. Here, the theory of categorical algebra, which is developed to internalize algebraic structures within a general category $\C$, comes to our aid; choosing the category of sets $\mathrm{Set}$ as $\C$ yields ordinary algebraic structures. However, $\C$ needs to have specific properties and cannot be a general category; in particular, in this paper, $\C$ is always $\mathrm{Cat}(\Ab)$. We will review some aspects of this theory in \Cref{Chrecall}. 

In particular, it is possible to define the concept of an $\underline{R}$-module within $\mathrm{Cat}(\Ab)$ and the category of $\underline{R}$-modules within $\mathrm{Cat}(\Ab)$. We denote this category by $\mathrm{LMod}_{\underline{R}}(\mathrm{Cat}(\Ab))$.
 In \Cref{defCategoryOfLeftModuleObjects}, we provide a mathematical definition of this category and refer to its objects not by the naive term $\underline{R}$-module within $\mathrm{Cat}(\Ab)$, but with the precise terminology of the category of left $\underline{R}$-module objects of $\mathrm{Cat}(\Ab)$. 

Roughly speaking, the first main result of this paper can be stated as follows:

\begin{displayquote}
\textit{the category of $R$-$\mathrm{Mod}$-enriched categories $\mathrm{Cat}(R$-$\mathrm{Mod})$ and the category of $\underline{R}$-modules within $\mathrm{Cat}(\Ab)$ $\mathrm{LMod}_{\underline{R}}(\mathrm{Cat}(\Ab))$ are equivalent.}
\end{displayquote}

See \Cref{thcategorialversionEnrichedModareMod} for the precise statement.

Thus, we transition from a module structure in the enrichment of categories to categories with an inherent module structure. 
This transition suggests that $R$-$\mathrm{Mod}$-enriched categories are $\mathrm{Cat}(\Ab)$-enriched functors, as it is well-known that modules exhibit this property.
We reiterate this fact starting with the ordinary algebraic concept of a module and then move to modules within $\mathrm{Cat}(\Ab)$. 

At the beginning of this introduction, we recalled that $\mathrm{Cat}$ is a $\mathrm{Cat}$-enriched category. Similarly, it is possible to define natural transformations between $\V$-enriched categories (see \Cref{defEnrichedNaturalTransformation}) and show that the category $\mathrm{Cat}(\V)$ is a $\mathrm{Cat}$-enriched category. We denote the hom-$\mathrm{Cat}$-object between two $\V$-enriched categories $\A$ and $\B$ by $\mathrm{Fun}^{\V}(\A,\B)$.

The following equivalence is well-known and holds almost tautologically:
\begin{equation}
    \label{eqAndreaBovo}
    \mathrm{Fun}^{\Ab}(\underline{R},\Ab) \simeq R\text{-}\mathrm{Mod};
\end{equation}
see \cite[Example 3.5.9]{riehl2014categorical}. This new perspective of $R$-modules as enriched functors is useful. For example, using \eqref{eqAndreaBovo}, the adjunctions of scalar extension-restriction and scalar restriction-coextension (see \cite[Example 7.6.9]{riehl2014categorical}) can be easily constructed.

Next, we consider modules within $\mathrm{Cat}(\Ab)$; that is, we aim to replace $\Ab$ with $\mathrm{Cat}(\Ab)$ in \eqref{eqAndreaBovo}. To do so, $\mathrm{Cat}(\Ab)$ must possess the same characteristics as $\Ab$. 
The hom-objects of an $\Ab$-enriched category have another canonical structure. Just as $\mathrm{Cat}$ is self-enriched (i.e., it is a $\mathrm{Cat}$-enriched category), $\mathrm{Cat}(\Ab)$ is a $\mathrm{Cat}(\Ab)$-enriched category, see \Cref{Chrecall}.
For the next result, we need a category of $\mathrm{Cat}(\Ab)$-enriched functors, so the source cannot be $\underline{R}$. It must be an internalization of a unitary ring $R$ in the category $\mathrm{Cat}(\mathrm{Cat}(\Ab))$.               
Using a similar approach as before, we can define the $\mathrm{Cat}(\Ab)$-enriched category $\2un{R}$, which has only one object $\{*\}$ and whose only non-trivial hom-object is $\Hom(*,*)=\underline{R}$.

The second main result can be stated as follows:

\begin{displayquote}
\textit{the category of $\mathrm{Cat}(\Ab)$-enriched functors $\mathrm{Fun}^{\mathrm{Cat}(\Ab)}(\2un{R},\mathrm{Cat}(\Ab))$ and the category of left $\underline{R}$-module objects of $\mathrm{Cat}(\Ab)$ $\mathrm{LMod}_{\underline{R}}(\mathrm{Cat}(\Ab))$ are equivalent.}
\end{displayquote}

See \Cref{thcategorialversionAction} for the precise statement.
This second result is the \enquote{with more objects} version of \eqref{eqAndreaBovo}.
This perspective of modules within a category as enriched functors is widely used in algebraic geometry and representation theory. For example, in Tannakian category theory, it is important that the category of linear representations of a group $G$ on finite-dimensional $k$-vector spaces, $\mathrm{Rep}(G)$, which is a category of enriched functors, is equivalent to the category of comodules, which is a category of modules within a suitable category; see \cite{DeligneTanna, deligne2009ii}.  


To summarize, the main results of this article can be stated as a single result as follows.
\begin{theorem}
\label{sumtheorem}
Let $R$ be a commutative and unitary ring. There is a chain of equivalences of categories
\[ Cat(R\mathrm{-}\mathrm{Mod})\simeq \mathrm{LMod}(Cat(\Ab))\simeq Fun^{Cat(\Ab)}(\2un{R},Cat(\Ab)).\]
\end{theorem}

\section*{Related Works}


This paper is the first in a series of articles \cite{DoniklinearMorita, DoniHigherCategorical} aimed at developing and comparing the Morita theory for dg-categories over $k$ and $k$-linear $\infty$-categories in an $\infty$-categorical setting. This article may initially seem unrelated to the purpose because it is written using category theory rather than $\infty$-category theory. However, it stems from the author's need to understand how being an $R$-linear category (or $R$-$\mathrm{Mod}$-enriched category) for a suitable unitary ring $R$ relates to being a $\Z$-linear category (or preadditive or $\Ab$-enriched category). 
Indeed, in \cite{DoniHigherCategorical}, we prove the $\infty$-categorical version of the main results of this paper, \Cref{thcategorialversionAction,thcategorialversionEnrichedModareMod}, which will be fundamental for achieving the final goal of this series in \cite{DoniklinearMorita}.
$\\$

I initially did not intend to write this article because I expected to find its results already documented somewhere. To my surprise, all my searches were inconclusive. To this day, I still believe the results are recorded somewhere or are part of mathematical folklore. I think my article can be useful because it can serve as an easy-to-find reference.

\section*{Structure of the Paper}
This article is designed to be as self-contained as possible and is split into two sections. In \Cref{Chrecall}, we define, recall, and construct all mathematical concepts and structures needed to state and prove the results. A reader familiar with enriched category theory and algebraic category theory may skip this part and refer to it only for clarifying notations, which are fairly standard. 

In \Cref{ChResults}, we prove the two main results of this work, \Cref{thcategorialversionAction,thcategorialversionEnrichedModareMod}, and discuss a possible generalization and a specific case.
The strategies of the proofs are purely categorical and classical in homological algebra. We will use \enquote{diagram chasing and pasting-diagram techniques}. The main difficulty in the proofs will be to demonstrate that such complex diagrams are commutative.

\section*{Acknowledgement}

I extend my gratitude to Paolo Stellari for his guidance throughout the recent years, which has enabled me to consolidate my ideas.
Furthermore, I express my appreciation to Sandra Mantovani, who, during a conversation, highlighted a crucial point in the proof of \Cref{thcategorialversionEnrichedModareMod}.

\section{Recalls}
\label{Chrecall}
\subsection{Notion of Enriched Category theory}
\label{chNotionofEnrichedCategoryTheory}
In this subsection, we define what is for us an enriched category, an enriched functor, and an enriched natural transformation: our definitions are inspired by \cite[\S 3]{riehl2014categorical}.

For this subsection, we fix a closed symmetric monoidal category $(\V,\otimes_{\V}, \I_{\V})$.

\begin{definition}
\label{defEnichedCategory}
A \textit{$\V$-enriched category $\D$} consists of
\begin{itemize}
    \item a collection of objects $x, y, z\in  \D$;
    \item for each pair $x, y \in \D$, a \textit{hom-object} $\D(x, y) \in V$;
    \item for each $x \in \D$, a morphism $id_x : \I_{\V}\to \D(x, x)$ in $\V$,
    \item for each triple $x, y, z \in \D$, a morphism $\circ:   \D(x, y)\otimes_{\V}\D(y, z) \to \D(x, z)$ in $\V$.
\end{itemize}
Such that the following diagrams commute for all $x, y, z, w \in \D$: 
\begin{equation}
\label{eqFirstDiagrammDefEnrichedCate}
\begin{tikzcd}[ampersand replacement=\&]
	{\D(x,y)\otimes\D(y,z)\otimes\D(z,w)} \&\&\& {\D(x,z)\otimes\D(z,w)} \\
	{\D(x,y)\otimes\D(y,w)} \&\&\& {\D(x,w)}
	\arrow["{ \circ\otimes id}"', from=1-1, to=1-4]
	\arrow["{id\otimes \circ}", from=1-1, to=2-1]
	\arrow["\circ", from=2-1, to=2-4]
	\arrow["\circ"', from=1-4, to=2-4]
\end{tikzcd}
\end{equation}

\begin{equation}
\label{eqSecondDiagrammDefEnrichedCate}
\begin{tikzcd}[ampersand replacement=\&]
	{\I_{\V}\otimes\D(y,w)} \& {\D(y,y)\otimes\D(y,w)} \\
	\& {\D(y,w)}
	\arrow["{id_{y}\otimes id}", from=1-1, to=1-2]
	\arrow["\circ"', from=1-2, to=2-2]
	\arrow["\cong"', from=1-1, to=2-2]
\end{tikzcd}
\end{equation}

\begin{equation}
\label{eqThirdDiagrammDefEnrichedCate}
\begin{tikzcd}[ampersand replacement=\&]
	{\D(y,w)\otimes\D(w,w)} \& {\D(y,w)\otimes\I_{\V}} \\
	{\D(y,w)}
	\arrow["{id\otimes id_{w}}"', from=1-2, to=1-1]
	\arrow["\cong", from=1-2, to=2-1]
	\arrow["\circ"', from=1-1, to=2-1]
\end{tikzcd}
\end{equation}

\end{definition}   

\begin{notation}
\label{notComposition}
Let $\D$ be a $\V$-enriched category and let $x,y$ be two objects of $\D$. We call \textit{the arrows (or morphism) of $\D$ from $x$ to $y$ } the elements of the set $\Hom_{\V}(\I_{\V},\D(x,y))$.
 Let $f\in \Hom_{\V}(\I_{\V},\D(x,y))$ be one such arrow and $z\in\D$ be an object, we call \textit{ precomposition with $f$} the arrow of $\V$:  
\[f^{*}:=\D(y,z)\cong \I_{\V}\otimes\D(y,z)\xrightarrow{f\otimes id} \D(x,y)\otimes\D(y,z)\xrightarrow{\circ}  \D(x,z).\]

We call \textit{ postcomposition with $f$} the arrow of $\V$: 
\[f_{*}:=\D(z,x)\cong \D(z,x)\otimes\I_{\V}\xrightarrow{id\otimes f} \D(z,x)\otimes\D(y,x) \xrightarrow{\circ}  \D(z,x).\]    
\end{notation}

\begin{definition}
\label{defEnrichedFunctor}
A \textit{$\V$-enriched functor} $F : \C\to \D$ between $\V$-enriched categories is given by an object
map \[ob(F):ob(\C)\to ob(\D):c\mapsto F(c)\] together with arrow
\[\C(x, y)\xrightarrow{F_{x,y}}\D(Fx, Fy)\]
in $\V$ for each $x, y \in \C$ such that the following diagrams commute for all $x, y, z \in \C$:

\begin{equation}
\label{eqFirstDigramdefVenrichedFunctor}
\begin{tikzcd}[ampersand replacement=\&]
	{\C(x,y)\otimes\C(y,z)} \&\& {\C(x,z)} \\
	{\D(Fx,Fy)\otimes\D(Fy,Fz)} \&\& {\D(Fx,Fz)}
	\arrow["\circ"', from=1-1, to=1-3]
	\arrow["{F_{x,y}\otimes F_{y,z}}", from=1-1, to=2-1]
	\arrow["\circ", from=2-1, to=2-3]
	\arrow["{F_{x,z}}"', from=1-3, to=2-3]
\end{tikzcd}
\end{equation}

\begin{equation}
\label{eqSecondDigramdefVenrichedFunctor}
\begin{tikzcd}[ampersand replacement=\&]
	{\I_{\V}} \&\& {\C(x,x)} \\
	\&\& {\D(Fx,Fx)}
	\arrow["{id_{x}}"', from=1-1, to=1-3]
	\arrow["{F_{x,x}}"', from=1-3, to=2-3]
	\arrow["{id_{Fx}}"', from=1-1, to=2-3]
\end{tikzcd}
\end{equation}
The commutativity of the last diagrams is sometimes called \textit{the functoriality axioms.} 
\end{definition}

The $\V$-enriched categories and the $\V$-enriched functors form a category, see \cite{kelly1982basic}.

\begin{notation}
\label{notCatV}
We denote by $Cat(\V)$ the category whose objects are $\V$-enriched categories and whose morphisms are the $\V$-enriched functors.
As usual, we call $Cat(\V)$ the \textit{the category of $\V$-enriched categories.}
\end{notation}

\begin{definition}
\label{defEnrichedNaturalTransformation}
 A \textit{$\V$-enriched natural transformation} $\alpha: F \Rightarrow G$ between a pair of $\V$-enriched
functors $F,G: C \rightrightarrows D$ consists of a morphism $\alpha_x : \I_{\V} \to \D(Fx,Gx)$ in $\V$ for each $x \in C$ such that for all $x, y \in C$ the following diagram commutes:
\[\begin{tikzcd}[ampersand replacement=\&]
	{\C(x,y)} \&\& {\D(Fx,Fy)} \\
	\\
	{\D(Gx,Gy)} \&\& {\D(Fx,Gy),}
	\arrow["{F_{x,y}}", from=1-1, to=1-3]
	\arrow["{G_{x,y}}"', from=1-1, to=3-1]
	\arrow["{(\alpha_{x})^{*}}"', from=3-1, to=3-3]
	\arrow["{(\alpha_{y})_{*}}", from=1-3, to=3-3]
\end{tikzcd}\]
Here, the starred arrows are precompositions or postcompositions.
\end{definition}

\subsection{Notions of Categorical Algebra}
In this short subsection, we recall some definitions from categorical algebra.

\begin{definition}[Left $\V$-left tensor category]
\label{defLeftTensorCategoriordinal}
Let $(\V,\otimes_{\V}, \I_{\V})$ be a monoidal category.
A $\V$-left tensored category is a pair $(\W, \otimes)$, where $\W$ is a category and $\otimes:\V\times \W\to \W$ is a bifunctor
such that the following axioms are satisfied: 

\begin{enumerate}
    \item The following diagram commutes:
\[\begin{tikzcd}[ampersand replacement=\&]
	\V\times\V\times\W \&\& {\V\times \W} \\
	\\
	\V\times\W \&\& {\W;}
	\arrow["{\otimes_{\V}\times id_{W}}", from=1-1, to=3-1]
	\arrow["\otimes", from=3-1, to=3-3]
	\arrow["{id\times \otimes}", from=1-1, to=1-3]
	\arrow["\otimes"', from=1-3, to=3-3]
\end{tikzcd}\]
\item the composition \[\W \simeq \mathbf{1} \times \W \xrightarrow{I_{\V}\times id_{\W}} \V\times \W   \xrightarrow{\otimes} \W\] is the identity.

\end{enumerate}

\end{definition}

We further need the following definitions:

\begin{definition}
\label{defMonoids}
Let $(\V,\otimes_{\V}, \I_{\V})$ be a monoidal category. A \textit{monoid in $\V$} is a triple $(w,m_{w}:w\otimes w\to w,u_{w}:\I_{\V}\to w)$ whose the first component is an object of $\V$ and the remaining components are arrows in $\V$ such that the following diagrams commute:
\[\begin{tikzcd}[ampersand replacement=\&]
	{w\otimes w\otimes w} \&\& {w\otimes w} \\
	\\
	{w\otimes w} \&\& {w;} \\
	{w\simeq w\otimes \I_{\V}} \& {w\otimes w} \& {w\simeq \I_{\V}\otimes w} \\
	\& {w.}
	\arrow["{id\otimes m_{w}}", from=1-1, to=3-1]
	\arrow["{m_{v}}", from=3-1, to=3-3]
	\arrow["{m_{w}\otimes id}"', from=1-1, to=1-3]
	\arrow["{m_{w}}"', from=1-3, to=3-3]
	\arrow["{id\otimes I_{w}}", from=4-1, to=4-2]
	\arrow["{I_{w}\otimes id}"', from=4-3, to=4-2]
	\arrow["id"', from=4-1, to=5-2]
	\arrow["id", from=4-3, to=5-2]
	\arrow["{m_{w}}", from=4-2, to=5-2]
\end{tikzcd}\]
As usual, we call $m_{w}$ \textit{the multiplication morphism of $w$} and $u_{w}$ the \textit{unit morphism of $w$}.

Furthermore, we denote by $Mon(\V)$ the category whose objects are monoids in $\V$ and an arrow $\hat{f}:(w,m_{w},u_{w})\to (v,m_{v},u_{v})$ of $Mon(\V)$ is an arrow $f:w\to v$ of $\V$ such that the following diagrams are commutative:

\begin{equation}
\label{eqDiagramArrowMonoids}
\begin{tikzcd}[ampersand replacement=\&]
	{w\otimes w} \&\& w \&\& {\I_{\V}} \& w \\
	\&\&\&\&\& {v.} \\
	{v\otimes v} \&\& {v;}
	\arrow["{f\otimes f}", from=1-1, to=3-1]
	\arrow["{m_{v}}", from=3-1, to=3-3]
	\arrow["{m_{w}}"', from=1-1, to=1-3]
	\arrow["f"', from=1-3, to=3-3]
	\arrow["{u_{w}}", from=1-5, to=1-6]
	\arrow["{u_{v}}"', from=1-5, to=2-6]
	\arrow["f", from=1-6, to=2-6]
\end{tikzcd}
\end{equation}
We call $f$ the \textit{underlying arrow of $\hat{f}$}.    
\end{definition}

\begin{definition}
\label{defCommutativeMonoid}
Let $(\V,\otimes_{\V}, \I_{\V})$ be a symmetric monoidal category.
We say that a monoid $w$ in $\V$ is a \textit{commutative monoid} if the following diagram is commutative:
\[\begin{tikzcd}[ampersand replacement=\&]
	{w\otimes w} \&\& {w\otimes w} \\
	\& {w;}
	\arrow["{\tau \cong}", from=1-1, to=1-3]
	\arrow["{m_{w}}", from=1-3, to=2-2]
	\arrow["{m_{w}}"', from=1-1, to=2-2]
\end{tikzcd}\]
here, $\tau$ is the braiding of the symmetric monoidal category $\V$.

Furthermore, we denote by $CMon(\V)\subseteq Mon(\V)$ the full subcategory spanned by commutative monoids in $\V$.   
\end{definition}

\begin{definition}
\label{defCategoryOfLeftModuleObjects}
Let $\V$ be a symmetric monoidal category.
Let $(z,m_{z}:z\otimes z\to z,u_{z}:\I_{\V}\to z)$ be a commutative monoid in $\V$ and let $(\W,\otimes)$ be a $\V$-left tensored category. 
We denote by $\mathrm{LMod}_{z}(\W)$ the category whoses objects are pairs $(w,m_{w}:z\otimes w\to w)$ whose first component is an object of $\W$ and whose second component is an arrow of $\W$ with source $z\otimes w$ and target $w$, satisfying following axioms:
\begin{itemize}
    \item[(1)] the following diagram is commutative:
\[\begin{tikzcd}[ampersand replacement=\&]
	z\otimes z\otimes w \&\& {z\otimes w} \\
	\\
	z\otimes w \&\& {w;}
	\arrow["{m_{z}\otimes id_{w}}", from=1-1, to=3-1]
	\arrow["m_{w}", from=3-1, to=3-3]
	\arrow["{id\otimes m_{w}}", from=1-1, to=1-3]
	\arrow["m_{w}"', from=1-3, to=3-3]
\end{tikzcd}\]
\item[(2)] the composition \[ w \simeq \I_{V} \otimes w \xrightarrow{u_{z}\otimes id_{w}} z\otimes w   \xrightarrow{m_{w}} w;\] is the identity.
\end{itemize}
An arrow of $\mathrm{LMod}_{z}(\W)$ is a pair $(f,\square)$ where the first component is an arrow of $\W$ and the second component is a commutative square:

\[\begin{tikzcd}[ampersand replacement=\&]
	{z\otimes w} \&\& w \\
	\\
	{z\otimes v} \&\& {v;}
	\arrow["{z\otimes f}", from=1-1, to=3-1]
	\arrow["{m_{v}}", from=3-1, to=3-3]
	\arrow["{m_{w}}"', from=1-1, to=1-3]
	\arrow["f"', from=1-3, to=3-3]
\end{tikzcd}\]
such that the following diagrams commute:
\begin{equation}
\label{eqMorphismOfModuleassociativity}
\begin{tikzcd}[ampersand replacement=\&]
	{z\otimes z\otimes w} \&\& {z\otimes w} \\
	\& {z\otimes z\otimes v} \&\& {z\otimes v} \\
	{z\otimes w} \&\& w \\
	\& {z\otimes v} \&\& {v.}
	\arrow[from=1-1, to=1-3]
	\arrow[from=1-1, to=3-1]
	\arrow[from=3-1, to=3-3]
	\arrow[from=1-3, to=3-3]
	\arrow[from=1-1, to=2-2]
	\arrow[from=1-3, to=2-4]
	\arrow[from=2-4, to=4-4]
	\arrow[from=2-2, to=4-2]
	\arrow[from=3-1, to=4-2]
	\arrow[from=4-2, to=4-4]
	\arrow[from=2-2, to=2-4]
	\arrow[from=3-3, to=4-4]
\end{tikzcd}
\end{equation}

\begin{equation}
\begin{tikzcd}[ampersand replacement=\&]
	{w\cong \I_{\V}\otimes w} \& {z\otimes w} \& w \\
	{v\cong \I_{\V}\otimes v} \& {z\otimes v} \& {v.}
	\arrow["{m_{w}}"', from=1-2, to=1-3]
	\arrow["{id\otimes v}", from=2-1, to=2-2]
	\arrow["{m_{v}}", from=2-2, to=2-3]
	\arrow["f"{description}, from=1-3, to=2-3]
	\arrow["{z\otimes f}"', from=1-2, to=2-2]
	\arrow["{id\otimes w}"', from=1-1, to=1-2]
	\arrow["f", from=1-1, to=2-1]
	\arrow["id"', curve={height=-18pt}, from=1-1, to=1-3]
	\arrow["id", curve={height=18pt}, from=2-1, to=2-3]
\end{tikzcd}
\end{equation}
We call $\mathrm{LMod}_{z}(\W)$ the \textit{category of left $z$-module objects of $\W$}.
The arrows of $\mathrm{LMod}_{z}(\W)$ are usually called \textit{$z$-equivariant arrows}.
\end{definition}

\begin{notation}
In the above definition, if $z$ is a commutative monoid, we sometimes denote the category $\mathrm{LMod}_z(\mathcal{\W})$ as $\mathrm{Mod}_z(\mathcal{W})$. Moreover, we call it the category of $z$-module objects of $\W$
\end{notation}

\begin{example}
Let $\Ab$ be the monoidal category of abelian groups.
A (commutative) monoid in $\Ab$ is a (commutative) ring with unit; let $R$ be one of them.
The category of left $R$-module objects of $\Ab$, $\mathrm{LMod}_{R}(\Ab)$, is equivalent to the category of $R$-$\mathrm{Mod}$. 
In this paper, the above two notions are used interchangeably.
\end{example}

\begin{remark}
It is possible to make the above definition more general and let $z$ vary.   
\end{remark}

\subsection{Structures and construction}

This subsection recalls and defines some mathematical structures and constructions related to the previous mathematical notions.

For this subsection, we fix a closed symmetric monoidal category $(\V,\otimes_{\V}, \I_{\V})$.

The small $\V$-enriched categories, the $\V$-enriched functors, and the $\V$-natural transformations assemble into a $2$-category (see \cite[\S 3.5]{riehl2014categorical} or \cite{kelly1982basic} for a detailed definition of this $2$-category).

\begin{notation}
\label{notEnrichedFunctor}
We denote by $Fun^{Cat(\V)}(\C,\D)$ the corresponding hom-$Cat$-object for each pair of $\V$-enriched categories $\C$ and $\D$.
\end{notation}

Moreover, the category of $\V$-enriched categories $Cat(\V)$ has a closed symmetric monoidal structure; we denote by $-\otimes^{\V}-$ its tensor and by $(-)^{(-)}$ its internal-hom. Let $\C$ and $\D$ be two $\V$-enriched categories, the tensor product between $\C$ and $\D$ is the $\V$-category $\C\otimes^{\V}\D$ whose objects are pairs $(c,d)\in \C\times\D$ and, for each pair  of objects $(c,d)$ and $(c',d')$ in $\C\otimes^{\V}\D$,  the corresponding hom-$\V$-object is 
\[C\otimes^{\V} \D((c, d); (c', d')):= C(c, c')\otimes_{\V}  D(d, d').\]
Modulo size issues, the internal-hom between two $\V$-enriched small categories $\C$ and $\D$ is the $\V$-enriched category $\D^{\C}$ whose objects are $\V$-enriched functors $F;G: \C \rightrightarrows \D$ and whose hom-objects are defined by the formula
\begin{equation}
\label{eqinternalhomcategory}
\D^{\C}(F,G) = \int_{c\in\C} \underline{D}(Fc,Gc), 
\end{equation}
where $\int$ is the enriched end, see \cite[\S 7.3]{riehl2014categorical}.

\begin{example}
    If $\V=Set$ we obtain the usual closed symmetric monoidal structure of $Cat=Cat(Set)$.
\end{example}

For the sake of completeness, we recall a specific case of Leibniz construction: for a detailed presentation see \cite[\S C.2]{RiehlElements}. Our definition instead is a specific case of \cite[C.2.8]{RiehlElements}.

Before we need to define the category of arrows of a small category.
\begin{notation}
    We denote with $\mathbb{2}$ the category whose the set of objects is $\{0,1\}$ and with only one non-trivial arrow with source $0$ and target $1$, $\{0\to 1\}$. 
\end{notation}

\begin{definition}
Let $\C$ be a category. We call \textit{category of arrows of $\C$} the category $\C^{\mathbb{2}}$.
\end{definition}

\begin{definition}[{Leibniz construction}]
\label{defLeibnizConstruction}
Let $(\mathcal{M}, \otimes, \I_{\mathcal{M}},{(-)}^{(-)})$ be a closed symmetric category admitting pushout and pullback, whose second component is the tensor bifunctor:
\[-\otimes-:\mathcal{M} \times \mathcal{M}\to \mathcal{M}; \] 
and, the last component is the internal-hom bifunctor:
\[{(-)}^{(-)}:\mathcal{M}^{op}\times\mathcal{M}\to \mathcal{M}. \]
There are two bifunctor:  
\[\hat{-\otimes -}:\mathcal{M}^{\mathbb{2}} \times \mathcal{M}^{\mathbb{2}}\to \mathcal{M}^{\mathbb{2}}\; , \hat{{(-)}^{(-)}}:(\mathcal{M}^{\mathbb{2}})^{op}\times \mathcal{M}^{\mathbb{2}}\to \mathcal{M}^{\mathbb{2}}; \] between the arrows category of $\mathcal{M}$. The bifunctor $\hat{-\otimes-}$ is called \textit{Leibniz tensor} and it carries a pair of objects $i:m\to m'\in \mathcal{M}^{2}$ and $j:n\to n'\in\mathcal{M}^{\mathbb{2}}$ in the objects $i\hat{\otimes}j\in \mathcal{M}^{\mathbb{2}}$ induced by the pushout diagram:
\[\begin{tikzcd}[ampersand replacement=\&]
	{m\otimes n} \& {m'\otimes n} \\
	{m\otimes n'} \& \bullet \\
	\&\& {m'\otimes n'.}
	\arrow["{m\otimes j}"', from=1-1, to=2-1]
	\arrow[from=2-1, to=2-2]
	\arrow["{i\otimes n}", from=1-1, to=1-2]
	\arrow[from=1-2, to=2-2]
	\arrow["{i\hat{\otimes}j}"', from=2-2, to=3-3]
	\arrow["{i\otimes n'}"', curve={height=12pt}, from=2-1, to=3-3]
	\arrow["{m'\otimes j}"', curve={height=-18pt}, from=1-2, to=3-3]
	\arrow["\lrcorner"{anchor=center, pos=0.125, rotate=180}, draw=none, from=2-2, to=1-1]
\end{tikzcd}\]
The bifunctors $\hat{{(-)}^{(-)}}$ is called \textit{Leibniz internal-hom} and it carries a pair of objects $i:m\to m'\in \mathcal{M}^{2}$ and $j:n\to n'\in\mathcal{M}^{\mathbb{2}}$ in the objects $\hat{m^k}\in \mathcal{M}^{\mathbb{2}}$ induced by the pushout diagram:
\[\begin{tikzcd}[ampersand replacement=\&]
	{n^{m'}} \\
	\& \bullet \& {n^{m}} \\
	\& {{n'}^{m'}} \& {{n'}^{m}.}
	\arrow["{\hat{j^i}}", from=1-1, to=2-2]
	\arrow[from=2-2, to=3-2]
	\arrow["{{n'}^{i}}", from=3-2, to=3-3]
	\arrow[from=2-2, to=2-3]
	\arrow["{j^m}"', from=2-3, to=3-3]
	\arrow["{j^{m'}}", curve={height=18pt}, from=1-1, to=3-2]
	\arrow["\lrcorner"{anchor=center, pos=0.125}, draw=none, from=2-2, to=3-3]
	\arrow["{n^i}"', curve={height=-18pt}, from=1-1, to=2-3]
\end{tikzcd}\]
\end{definition}

\begin{example}
$Cat(\Ab)$ satisfies the hypothesis of the above construction. Indeed we already know that it has a closed symmetric structure. In \cite{kelly2001v}, \cite{wolff1974v} or \cite{MuroHomOp} the authors prove that it is cocomplete; in particular, it admits pushout.
Moreover, as a consequence of \cite[\S ]{MuroHomOp} or \cite[Proposition 5.6]{borceux1994handbook}, it is also complete; in particular, it admits pullback. 
\end{example}

\begin{notation}
\label{notCategorywithOnlyOneObject}
Let $M$ be a monoid in $\V$.
We denote by $\underline{M}$ the $\V$-enriched category with only one object $*$ and whose the only non-trivial hom-$\V$-object is $\underline{M}(*,*)=M$ and the only non-trivial composition morphism of $\un{M}$ is the multiplication of $M$
\[m_{M}:Hom(*,*)\otimes_{\V} Hom(*,*)= M\otimes_{\V} M\to Hom(*,*)=M. \]
\end{notation}

\begin{definition}
\label{defUnderlineFunctorOrdinal}
We call \textit{underline functor} both of the following canonical functors: 
\[\underline{(-)}: Mon(\V)\to Cat(\V): M\mapsto\underline{M};\]
\[\underline{(-)}: CMon(\V)\hookrightarrow Mon(\V)\to Cat(\V): M\mapsto\underline{M};\]
here, $\underline{M}$ is as in \Cref{notCategorywithOnlyOneObject}.
\end{definition}

Since the monoidal structure in $\mathcal{V}$ is symmetric, the categories $Mon(\mathcal{V})$ and $CMon(\mathcal{V})$ have a symmetric monoidal structure with the property that the forgetful functors 
\[CMon(\mathcal{V}) \to \mathcal{V};\]
and,
\[Mon(\mathcal{V}) \to \mathcal{V};\] are monoidal. In both the structures the unit is the unit of $\mathcal{V}$ with the canonical multiplication, and the tensor functors are defined as follows:
\begin{equation}
\begin{split}
-\otimes- : Mon(\mathcal{V}) \times Mon(\mathcal{V}) & \to Mon(\mathcal{V}) \\
((A, m), (B, n)) & \mapsto (A, m) \otimes (B, n) = (A \otimes_{\mathcal{V}} B, mn);
\end{split}
\end{equation}
and
\begin{equation}
\begin{split}
-\otimes -: CMon(\mathcal{V}) \times CMon(\mathcal{V}) & \to CMon(\mathcal{V}) \\
((A, m), (B, n)) & \mapsto (A, m) \otimes (B, n) = (A \otimes_{\mathcal{V}} B, mn);
\end{split}
\end{equation}
where the second components $mn$ are both given by:
\[ mn : (A \otimes_{\mathcal{V}} B) \otimes_{\mathcal{V}} (A \otimes_{\mathcal{V}} B) \cong (A \otimes_{\mathcal{V}} A) \otimes_{\mathcal{V}} (B \otimes_{\mathcal{V}} B) \xrightarrow{m \otimes_{\mathcal{V}} n} A \otimes_{\mathcal{V}} B. \]
Note that in the definition of $mn$, it is important that the monoidal category $\mathcal{V}$ is symmetric.

For a detailed construction of the above symmetric monoidal structure see \cite{porst2008categories} or \cite{joyal1993braided}. These monoidal structures are studied also in \cite{johnstone2002sketches}.

For each symmetric monoidal category $\mathcal{V}$, the underline functor is monoidal. In fact, using the definition of the monoidal structure of $Cat(\mathcal{V})$, the following equivalence holds:
\[ \underline{M} \otimes^{\mathcal{V}} \underline{M} \simeq \underline{M \otimes_{\mathcal{V}} M}. \]
In particular, it implies that $\underline{(-)}$ sends monoids into monoids, so it induces a functor:
\[Mon(\underline{(-)}) : Mon(Mon(\V))\to Mon(Cat(\V)).\]

There is a canonical inclusion
\begin{equation}
\label{eqJInclusion}
\iota: Mon(\V)\to Mon(Mon(\V)),
\end{equation} which carries an object of $Mon(\V)$ represented by the pair $(w,m_{w}:w\otimes w\to w)$ into the object of $Mon(Mon(\V))$ represented by the pair $((w,m_{w}:w\otimes w\to w), \square_{w})$ where $\square_{w}$ is the arrow in $Mon(\V)$ designed in the next diagram:
\[\begin{tikzcd}[ampersand replacement=\&]
	{(w\otimes w)\otimes (w\otimes w)} \& {w\otimes w} \\
	{w\otimes w} \& {w.}
	\arrow["{m_{w}\otimes m_{w}}", from=1-1, to=2-1]
	\arrow["{m_{w}\otimes m_{w}}"', from=1-1, to=1-2]
	\arrow["{m_{w}}"', from=1-2, to=2-2]
	\arrow["{m_{w}}", from=2-1, to=2-2]
\end{tikzcd}\]
Moreover, $\iota$ sends an arrow $\hat{f}:(w,m_{w},u_{w})\to (v,m_{v},u_{v})$ with $f:w\to v$ as the underlying arrow into the arrow $\iota(\hat{f})$ with $\hat{f}:(w,m_{w},u_{w})\to (v,m_{v},u_{v})$ as underlying arrow. $\iota(\hat{f})$ is an arrow of $Mon(Mon(\V))$ because the following diagrams are canonically commutative: 
\begin{equation}
\begin{tikzcd}[ampersand replacement=\&]
	\& {(w\otimes w)\otimes (w\otimes w)} \& {w\otimes w} \\
	{(v\otimes v)\otimes (v\otimes v)} \& {w\otimes w} \& w \\
	{v\otimes v} \& {v\otimes v} \\
	\& {v;}
	\arrow[from=1-2, to=2-2]
	\arrow[from=1-2, to=1-3]
	\arrow[from=1-3, to=2-3]
	\arrow[from=2-2, to=2-3]
	\arrow[from=1-2, to=2-1]
	\arrow[from=1-3, to=3-2]
	\arrow[from=2-2, to=3-1]
	\arrow[from=2-3, to=4-2]
	\arrow[from=3-1, to=4-2]
	\arrow[from=2-1, to=3-2]
	\arrow[from=3-2, to=4-2]
	\arrow[from=2-1, to=3-1]
\end{tikzcd}
\end{equation}

\begin{equation}
\begin{tikzcd}[ampersand replacement=\&]
	\& {w\otimes w} \& w \\
	\& {v\otimes v} \& {v.} \\
	{\I_{\V}\otimes\I_{\V}} \& {\I_{\V}}
	\arrow[from=1-2, to=1-3]
	\arrow[from=3-1, to=3-2]
	\arrow[from=1-3, to=2-3]
	\arrow[from=1-2, to=2-2]
	\arrow[from=3-1, to=2-2]
	\arrow[from=3-2, to=2-3]
	\arrow[from=2-2, to=2-3]
	\arrow[from=3-2, to=1-3]
	\arrow[from=3-1, to=1-2]
\end{tikzcd}
\end{equation}
By composing $\iota$ and $Mon(\underline{(-)})$ we obtain an arrow
\[Mon(\underline{(-)})\circ \iota :Mon(\V)\to Mon(Cat(\V)).      \] 

Since the monoidal structure of $Cat(\V)$ is symmetric, it has its underline functor.

\begin{definition}
\label{defDoubleUnderlineFuntorOrdinal}
We call \textit{the double-underline functor} both the functors:    
\begin{equation}
\label{eqOrdinalDuobleMon}
\2un{(-)}:=Mon(\V)\xrightarrow{Mon(\underline{(-)})\circ \iota} Mon(Cat(\V))\xrightarrow{\underline{(-)}} Cat(Cat(\V));  
\end{equation}
and its restriction:
\begin{equation}
\label{eqordinalDuobleCommutativo}
\2un{(-)}:=CMon(\V)\xrightarrow{Mon(\underline{(-)})\circ \iota} Mon(Cat(\V))\xrightarrow{\underline{(-)}} Cat(Cat(\V)).  
\end{equation}

\end{definition}

Now we specialize to the case when $(\V,\otimes_{\V},\I_{\V})=(\Ab,\otimes_{\Z}, \Z)$; this is the only case we are interested in. At the end of this paper, we spend some words on a possible generalization, see \Cref{rmkGeneralizzation}.

The category $Cat(\Ab)$ has a symmetric monoidal closed structure and the internal-hom obeys the formula (\ref{eqinternalhomcategory}).

It is known that (commutative) monoids of $\Ab$ are actually (commutative) rings, $Mon(\Ab)\simeq Ring $ ($CMon(\Ab)\simeq CRing $), 
thus the underline functors become:
\[\underline{(-)}: Ring\to Cat(\Ab): M\mapsto\underline{M};\]
and,
\[\underline{(-)}: CRing\to Cat(\Ab): M\mapsto\underline{M}.\]
While the double-underline functors become:
\[\2un{(-)}:Ring \to Cat(Cat(\Ab));\]
and
\[\2un{(-)}:CRing \to Cat(Cat(\Ab)).\]

We choose a commutative unitary ring  $R\in CRing$. 
In particular, for what we have observed above, $\underline{R}$ is a monoid in $Cat(\Ab)$ and $\2un{R}$ is a $Cat(\Ab)$-enriched category.

It is well-known that the triple $(R\text{-}\mathrm{Mod}, \otimes_{R}, R )$ forms a symmetric monoidal category, so the category $Cat(R\text{-}\mathrm{Mod})$ is well-defined. 

Now we describe the functor when we consider the category of $R$-modules as a category of left $R$-module objects.
Let $(\A,m_{\A})$ and $(\B,m_{\B})$ be two objects of $\mathrm{Mod}_{R}(\Ab)$, the module objects $\A\otimes_{R}\B$ is a pair $(\A\otimes_{R}\B,m_{\A\otimes_{R}\B})$ whose the first component is given by the following coequalizer diagram:
\begin{equation}
\label{eqCorrezioneErrore5}
\begin{tikzcd}[ampersand replacement=\&]
	{\A\otimes_{\Z}R\otimes_{\Z}\B} \& {\A\otimes_{\Z}\B} \& {coeq(m_{\A}\otimes id,id\otimes m_{\B})=\A\otimes_{R}\B;}
	\arrow["{id\otimes m_{\B}}", shift left=2, from=1-1, to=1-2]
	\arrow["{m_{\A}\otimes id}"', shift right=2, from=1-1, to=1-2]
	\arrow["q",from=1-2, to=1-3]
\end{tikzcd}
\end{equation}
instead $m_{\A\otimes_{R}\B}$ is the arrow induced by the universal property of coequalizers in the below diagram:
\[\begin{tikzcd}[ampersand replacement=\&,sep=small]
	{R\otimes_{\Z} R\otimes_{\Z}\A\otimes_{\Z}\B} \\
	{R\otimes_{\Z}\A\otimes_{\Z}\B} \\
	{\A\otimes_{\Z}R\otimes_{\Z}\B} \&\& {R\otimes_{\Z}(\A\otimes_{R}\B)} \& {coeq(R\otimes_{\Z}(m_{\A}\otimes_{\Z} id),R\otimes_{\Z}(id\otimes_{\Z} m_{\B}))} \\
	{\A\otimes_{\Z}\B} \&\& {\A\otimes_{R}\B.}
	\arrow["{id\otimes m_{\B}}", shift left=2, from=3-1, to=4-1]
	\arrow["{m_{\A}\otimes id}"', shift right=2, from=3-1, to=4-1]
	\arrow["q"', from=4-1, to=4-3]
	\arrow[from=3-1, to=3-3]
	\arrow["{m_{\A\otimes_{R}\B}}", dashed, from=3-3, to=4-3]
	\arrow["\simeq"{marking, allow upside down}, draw=none, from=3-1, to=2-1]
	\arrow["{R\otimes q}"', from=2-1, to=3-3]
	\arrow["{R\otimes_{\Z}(id\otimes m_{\B})}", shift left=3, from=1-1, to=2-1]
	\arrow["{R\otimes(m_{\A}\otimes id)}"', shift right=3, from=1-1, to=2-1]
	\arrow["\simeq"{marking, allow upside down}, draw=none, from=3-3, to=3-4]
\end{tikzcd}\]
It is crucial to note that to utilize the universal property of the coequalizers, it is imperative for the functor $R\otimes_{\Z}\text{-}$ to commute with small colimits. Additionally, in the provided definition of the tensor product between $R$-modules, the ring $R$ needs to be commutative. 

We have now established all the necessary mathematical concepts to state the results.


\section{Results}
\label{ChResults}
In this section, we prove the desired results, which are explained in the introduction.

\begin{theorem}
\label{thcategorialversionEnrichedModareMod}
Let $R$ be a commutative ring with unit. There is an equivalence of categories
\begin{equation*}
\mathrm{LMod}_{\underline{R}}(Cat(\Ab))\simeq Cat(R\text{-}\mathrm{Mod}).
\end{equation*}
\end{theorem}

\begin{proof}
Let us start by defining a functor:
\[
\phi: \mathrm{LMod}_{\underline{R}}(Cat(\Ab)) \to Cat^{R\text{-}\mathrm{Mod}}.
\]

Let $\mathcal{A} = (\mathcal{A}, m_{\mathcal{A}}: \underline{R} \otimes^{\Ab} \mathcal{A} \to \mathcal{A})$ be an object of $\mathrm{LMod}_{\underline{R}}(Cat(\Ab))$, it satisfies the following axiom: the following two diagrams are commutative:

\begin{equation}
\label{eqCorrezioneErrore7}
\begin{tikzcd}[ampersand replacement=\&]
	{\un{R}\otimes^{\Ab}\un{R}\otimes^{\Ab}\A} \& {\un{R}\otimes^{\Ab}\A} \\
	{\un{R}\otimes^{\Ab}\A} \& {\A;}
	\arrow[from=1-1, to=2-1]
	\arrow[from=1-1, to=1-2]
	\arrow[from=1-2, to=2-2]
	\arrow[from=2-1, to=2-2]
\end{tikzcd}
\end{equation}

\begin{equation}
\label{eqCorrezioneErrore8}
\begin{tikzcd}[ampersand replacement=\&]
	\A \& {\un{\Z}\otimes^{\Ab}\A} \& {\un{R}\otimes^{\Ab}\A} \& {\A.}
	\arrow[from=1-2, to=1-3]
	\arrow[from=1-3, to=1-4]
	\arrow["\simeq"{description}, draw=none, from=1-1, to=1-2]
	\arrow["id", curve={height=18pt}, from=1-1, to=1-4]
\end{tikzcd}
\end{equation}

By definition, 
$m_{\mathcal{A}}$ is composed of a morphism of sets
\[
ob(m_{\mathcal{A}}): ob(\underline{R} \otimes^{\Ab} \mathcal{A}) \simeq ob(\mathcal{A}) \to ob(\mathcal{A}): (a) \mapsto m_{\mathcal{A}}(a),
\]
and, for each $a,b \in \mathcal{A}$, by a morphism of abelian groups:
\begin{equation}
\label{eqActiononHomObject}
m_{\mathcal{A} (a,b)}: \underline{R}(*,*) \otimes_{\mathbb{Z}} \mathcal{A}(a,b) \to \mathcal{A}(m_{\mathcal{A}}(a), m_{\mathcal{A}}(b)),
\end{equation}
plus two axioms, see Definition \ref{defCategoryOfLeftModuleObjects}.

The axiom $(2)$ in \Cref{defCategoryOfLeftModuleObjects} implies that the following morphism of sets 
\[
ob(\mathcal{A}) \simeq ob(\underline{Z} \otimes^{\Ab} \mathcal{A}) \to ob(\underline{R} \otimes^{\Ab} \mathcal{A}) \simeq ob(\mathcal{A}) \to ob(\mathcal{A}): (a) \mapsto m_{\mathcal{A}}(a)
\]
is the identity: that is, for each $a \in \mathcal{A}$, the equation $m_{\mathcal{A}}(a) = a$ holds.

In particular, \eqref{eqActiononHomObject} becomes:
\[
m_{\mathcal{A} (a,b)}: \underline{R}(*,*) \otimes_{\mathbb{Z}} \mathcal{A}(a,b) \to \mathcal{A}(a,b).
\]
The ($a,b$)-components of \eqref{eqCorrezioneErrore7} and \eqref{eqCorrezioneErrore8} become: 
\[\begin{tikzcd}[ampersand replacement=\&]
	{R\otimes_{\Z}R\otimes_{\Z}\A(a,b)} \& {R\otimes_{\Z}\A(a,b)} \\
	{R\otimes_{\Z}\A(a,b)} \& {\A(a,b);}
	\arrow[from=1-1, to=2-1]
	\arrow[from=1-1, to=1-2]
	\arrow[from=1-2, to=2-2]
	\arrow[from=2-1, to=2-2]
\end{tikzcd}\]

and 
\[\begin{tikzcd}[ampersand replacement=\&]
	{\A(a,b)} \& {\Z\otimes_{\Z}\A(a,b)} \& {R\otimes_{\Z}\A(a,b)} \& {\A(a,b).}
	\arrow[from=1-2, to=1-3]
	\arrow[from=1-3, to=1-4]
	\arrow["\simeq"{description}, draw=none, from=1-1, to=1-2]
	\arrow["id", curve={height=18pt}, from=1-1, to=1-4]
\end{tikzcd}\]

So, the starting $(\mathcal{A}, m_{\mathcal{A}}: \underline{R} \otimes^{\Ab} \mathcal{A} \to \mathcal{A})$ has hom-$R$-$\mathrm{Mod}$-objects.
Now, we prove that it is an $R$-$\mathrm{Mod}$-enriched category, Definition \ref{defEnichedCategory}.
Let $a$ be an object of $\mathcal{A}$, the identity morphism 
$id_{a}: R \to \mathcal{A}(a,a)$ ($id_{a}$ is a morphism of $R$-modules) is defined by the following commutative diagram:
\begin{equation}
\label{eqUnit}
\begin{tikzcd}[ampersand replacement=\&]
	{R\otimes_{\Z}R} \&\&\& R \\
	{R\otimes_{\Z}(R\otimes_{\Z}\Z)} \&\&\& {R\otimes_{\Z}\Z} \\
	{R\otimes_{\Z}(R\otimes_{\Z}\A(a,a))} \&\&\& {R\otimes_{\Z}\A(a,a)} \\
	{R\otimes_{\Z}\A(a,a)} \&\&\& {\A(a,a).}
	\arrow["{m_{R}}", from=1-1, to=1-4]
	\arrow["\simeq"', from=1-1, to=2-1]
	\arrow["\simeq", from=1-4, to=2-4]
	\arrow["{id_{R}\otimes(id_{R}\otimes id_{a})}"', from=2-1, to=3-1]
	\arrow["{id_{R}\otimes id_{a}}", from=2-4, to=3-4]
	\arrow["{m_{R}\otimes_{\Z}id_{\A(a,a)}}", from=3-1, to=3-4]
	\arrow["{m_{R}\otimes_{\Z}id_{\A(a,a)}}", from=3-1, to=4-1]
	\arrow["{m_{\A(a,a)}}"', from=3-4, to=4-4]
	\arrow["{m_{\A(a,a)}}", from=4-1, to=4-4]
\end{tikzcd}
\end{equation}

Instead, the composition is defined as follows.
We start with the solid diagram below:
\begin{equation}
\label{eqComposition}
\begin{tikzcd}[ampersand replacement=\&]
	{R\otimes_{\Z}\A(a,b)\otimes_{\Z}R\otimes_{\Z}\A(b,c)} \\
	{R\otimes_{\Z} (\A(a,b)\otimes_{\Z}\A(b,c))} \&\&\& {R\otimes_{\Z}\A(a,c)} \\
	{\A(a,b)\otimes_{\Z}\A(b,c)} \&\&\& {\A(a,c);} \\
	\&\& {R\otimes_{\Z}(\A(a,b)\otimes_{R}\A(b,c))} \\
	\&\& {\A(a,b)\otimes_{R}\A(b,c)}
	\arrow["{id\otimes_{\Z} m_{\A(b,c)}}", shift left=2, from=2-1, to=3-1]
	\arrow["{m_{\A(a,b)}\otimes_{\Z}id}"', shift right=2, from=2-1, to=3-1]
	\arrow["{R\otimes_{\Z}(id\otimes_{\Z} m_{\A(b,c=})}", shift left=2, from=1-1, to=2-1]
	\arrow["{R\otimes_{\Z}(m_{\A(a,b)}\otimes_{\Z}id)}"', shift right=2, from=1-1, to=2-1]
	\arrow["q"', from=3-1, to=5-3]
	\arrow["{R\otimes_{\Z}q}"', from=2-1, to=4-3]
	\arrow["\circ"', from=3-1, to=3-4]
	\arrow["{R\otimes_{\Z}\circ}"', from=2-1, to=2-4]
	\arrow["{m_{\A(a,c)}}", from=2-4, to=3-4]
	\arrow[dashed, from=4-3, to=2-4]
	\arrow[dashed, from=5-3, to=3-4]
	\arrow["{m_{\A(a,b)\otimes_{R}\A(b,c)}}", from=4-3, to=5-3]
\end{tikzcd}
\end{equation}
where $\circ$'s are the compositions of $\mathcal{A}$. 
We claim that the following equation holds $-\circ- (id\otimes m_{b,c})= -\circ- (m_{a,b}\otimes id)$.    
Assuming for the moment that the claim is true.  Since $R\otimes_{\Z}-$ commutes with colimits, as a consequence of the universal property of coequalizers, we obtain the two dotted arrows.

In particular, the rectangle with two dotted arrows in \eqref{eqComposition} will be the composition of $\mathcal{A}$ when we consider it as an $R$-$\mathrm{Mod}$-enriched category.

Now we prove the claim. Since $\mathcal{A}$ is a left $\un{R}$-module object for each $a,b,c\in \mathcal{A}$, the following diagram commutes: 
\begin{equation}
\label{eqCorrezioneErrore}
\begin{tikzcd}[ampersand replacement=\&,sep=small]
	{R\otimes_{\Z}\A(a,b)\otimes_{\Z}R\otimes\A(b,c)} \& {R\otimes_{\Z}R\otimes_{\Z}\A(a,b)\otimes_{\Z}\A(b,c)} \&\& {R\otimes_{\Z}\A(a,c)} \\
	{\A(a,b)\otimes\A(b,c)} \&\&\& {\A(a,c)};
	\arrow["\simeq", draw=none, from=1-1, to=1-2]
	\arrow["\circ", curve={height=-18pt}, from=1-1, to=1-4]
	\arrow["{m_{\A(a,b)}\otimes_{\Z} m_{\A(b,c)}}", from=1-1, to=2-1]
	\arrow["\circ", from=2-1, to=2-4]
	\arrow["{m_{\A(a,c)}}", from=1-4, to=2-4]
	\arrow["{m_R\otimes \circ }"', from=1-2, to=1-4]
\end{tikzcd}
\end{equation}

In \eqref{eqCorrezioneErrore}, $(m_{R}\times \circ)\simeq =\circ$ holds because the definition of composition in $\un{R}\otimes^{\Ab}\A$.

We precompose the diagram \eqref{eqCorrezioneErrore} with the morphism \[\A(a,b)\otimes_{\Z}R\otimes_{\Z}\A(b,c)\simeq \Z\otimes_{\Z}\A(a,b)\otimes_{\Z}R\otimes_{\Z}\A(b,c)\xrightarrow{u_{R}\otimes_{\Z} id} R\otimes_{\Z}\A(a,b)\otimes_{\Z}R\otimes_{\Z}\A(b,c);\]  
and we obtain the commutative diagram:

\begin{equation}
\label{eqCorrezioneErrore1}
\begin{tikzcd}[ampersand replacement=\&,sep=small]
	{\A(a,b)\otimes_{\Z}R\otimes_{\Z}\A(b,c)} \\
	{\Z\otimes_{\Z}\A(a,b)\otimes_{\Z}R\otimes_{\Z}\A(b,c)} \&\&\&\& {\A(x,y)\otimes_{\Z}\A(b,c)} \\
	\& {R\otimes_{\Z}\A(a,b)\otimes_{\Z}R\otimes_{\Z}\A(b,c)} \\
	\\
	{\A(a,b)\otimes_{\Z}\A(b,c)} \&\&\&\& {\A(a,c).}
	\arrow["{u_{R}\otimes_{\Z}id}", from=2-1, to=3-2]
	\arrow["{id\otimes_{\Z}m_{\A(b,c)}}"', from=2-1, to=5-1]
	\arrow["\circ"', from=5-1, to=5-5]
	\arrow["{m_{\A(a,b)}\otimes_{\Z}id}", from=2-1, to=2-5]
	\arrow["\circ"', from=2-5, to=5-5]
	\arrow[from=3-2, to=5-5]
	\arrow["{m_{\A(a,b)}\otimes_{\Z}m_{\A(b,c)}}", from=3-2, to=5-1]
	\arrow["{m_{\A(a,b)}\otimes_{\Z}m_{\A(b,c)}}"', from=3-2, to=2-5]
	\arrow["\simeq"{marking, allow upside down}, draw=none, from=1-1, to=2-1]
\end{tikzcd}
\end{equation}

The commutative of the perimeter of \eqref{eqCorrezioneErrore} implies the claim.
Note that to pass from \eqref{eqCorrezioneErrore} to \eqref{eqCorrezioneErrore1} it is important that $R$ is unitary.

It remains to prove that the three diagrams in Definition \ref{defEnichedCategory} are commutative.
In this case \eqref{eqFirstDiagrammDefEnrichedCate} becomes:
\[\begin{tikzcd}[ampersand replacement=\&,sep=tiny]
	{R\otimes_{\Z}(\A(a,b)\otimes_{R}\A(b,c)\otimes_{R}\A(c,d))} \& {R\otimes_{\Z}(\A(a,c)\otimes_{R}\A(c,d))} \\
	\& {\A(a,b)\otimes_{R}\A(b,c)\otimes_{R}\A(c,d)} \& {\A(a,c)\otimes_{R}\A(c,d)} \\
	{R\otimes_{\Z}(\A(a,b)\otimes_{R}\A(b,d))} \& {R\otimes_{\Z}(\A(a,d))} \\
	\& {\A(a,b)\otimes_{R}\A(b,d)} \& {\A(a,d).}
	\arrow[from=1-1, to=2-2]
	\arrow[from=1-1, to=3-1]
	\arrow[from=3-1, to=4-2]
	\arrow[curve={height=12pt}, from=2-2, to=4-2]
	\arrow[from=1-1, to=1-2]
	\arrow[from=3-1, to=3-2]
	\arrow[curve={height=-18pt}, from=1-2, to=3-2]
	\arrow[from=4-2, to=4-3]
	\arrow[from=2-2, to=2-3]
	\arrow[from=2-3, to=4-3]
	\arrow[from=1-2, to=2-3]
	\arrow[from=3-2, to=4-3]
\end{tikzcd}\]
In the above diagram, all the arrows but \[R\otimes_{\Z}\A(a,c)\to\A(a,c);\] are defined using the universal property of coequalizers, so by uniqueness the diagram is commutative. 

We only prove the commutativity of \eqref{eqSecondDiagrammDefEnrichedCate} of the remaining diagrams because the proof of \eqref{eqThirdDiagrammDefEnrichedCate} is symmetric.
Putting all the information together, we have to prove that the horizontal lines in the below commutative diagram are identities:  
\[\begin{tikzcd}[ampersand replacement=\&]
	{R\otimes_{\Z}(R\otimes_{R}\A(b,c))} \&\& {R\otimes_{\Z}(\A(b,b)\otimes_{R}\A(b,c))} \&\& {R\otimes_{\Z}\A(b,c)} \\
	{\A(b,c)\simeq R\otimes_{R}\A(b,c)} \&\& {\A(b,b)\otimes_{R}\A(b,c)} \&\& {\A(b,c).}
	\arrow[from=1-1, to=1-3]
	\arrow[from=1-3, to=1-5]
	\arrow[from=1-1, to=2-1]
	\arrow[from=2-1, to=2-3]
	\arrow[from=1-3, to=2-3]
	\arrow[from=2-3, to=2-5]
	\arrow["{m_{\A(b,c)}}", from=1-5, to=2-5]
\end{tikzcd}\]
Since the upper horizontal lines are the lower tensored with $R\otimes_{\Z}-$, it is sufficient to prove that the lower line is the identity morphism.
By construction, we have the following commutative diagram:
\[\begin{tikzcd}[ampersand replacement=\&]
	{\A(b,c)} \& {(\Z\otimes_{\Z}R)\otimes_{R}\A(b,c)} \&\& {\A(b,b)\otimes_{R}\A(b,c)} \\
	\& {\Z\otimes_{\Z}(R\otimes_{R}\A(b,c))} \& {\Z\otimes_{\Z}(\A(b,c))} \& {\A(b,b)\otimes_{\Z}\A(b,c)} \& {\A(b,c).}
	\arrow["\simeq"{description}, draw=none, from=1-1, to=1-2]
	\arrow[from=2-3, to=2-4]
	\arrow["q"', from=2-4, to=1-4]
	\arrow["\simeq"{marking, allow upside down}, draw=none, from=1-2, to=2-2]
	\arrow[from=2-2, to=2-3]
	\arrow[from=1-2, to=1-4]
	\arrow["\circ"', from=2-4, to=2-5]
	\arrow["\circ", from=1-4, to=2-5]
\end{tikzcd}\]
The upper lines are the lines whose composition we have to prove is the identity.
But the composition of the lower lines is the identity because $\A$ is a $\Ab$-enriched category.
$\\$

Therefore we define the morphism of sets:
\[
ob(\psi): ob(\mathrm{LMod}_{\underline{R}}(Cat(\Ab))) \to ob(Cat^{R\text{-}\mathrm{Mod}}): \mathcal{A} \mapsto \psi(\mathcal{A}):=\mathcal{A},
\]
as a kind of identity morphism.
For convenience, we sometimes write $\mathcal{A}$ and $\mathcal{B}$ instead of $\psi(\mathcal{A})$ and $\psi(\mathcal{B})$.

For each pair $\mathcal{A},\mathcal{B}$ of objects of $\mathrm{LMod}_{\underline{R}}(Cat(\Ab))$, we define the morphism between the sets of objects again as a kind of identity:
\[
F_{\mathcal{A},\mathcal{B}}: \Hom_{\mathrm{LMod}_{\underline{R}}(Cat(\Ab)))}(\mathcal{A},\mathcal{B}) \to \Hom_{Cat^{R\text{-}\mathrm{Mod}}}(\mathcal{A},\mathcal{B}).
\]

Indeed, let $(f,\square)$ be elements of $\text{Hom}_{\mathrm{LMod}_{\underline{R}}(Cat(\Ab))}(\mathcal{A},\mathcal{B})$ 
where $\square$ is the following commutative diagram:  
\begin{equation}
\label{eqDajeChecisiamo}    
\begin{tikzcd}[ampersand replacement=\&]
	{\un{R}\otimes^{\Ab} \A} \&\& \A \\
	\\
	{\un{R}\otimes^{\Ab} \B} \&\& \B.
	\arrow["{m_{\A}}"', from=1-1, to=1-3]
	\arrow["{\un{R}\otimes f}"', from=1-1, to=3-1]
	\arrow["{m_{\B}}"', from=3-1, to=3-3]
	\arrow["f"', from=1-3, to=3-3]
\end{tikzcd}
\end{equation}
We define $\psi(f):\mathcal{A}\to\mathcal{B}$ as the $R$-$\mathrm{Mod}$-enriched functor which behaves on objects as $f$:
\[
ob(\psi(f))=ob(f): ob(\mathcal{A}) \to ob(\mathcal{B}): a \mapsto \psi(f)(a) = f(a);
\]
and for each pair $a,b\in\A$ the morphism $\psi_{a,b}$ in $R$-$\mathrm{Mod}$ is the ($a,b$)-component of \eqref{eqDajeChecisiamo}:

\begin{equation}
\label{eqCorrezioneErrore6}
\begin{tikzcd}[ampersand replacement=\&]
	{R\otimes_{\Z} \A(a,b)} \&\& {\A(a,b)} \\
	\\
	{R\otimes_{\Z} \B(f(a),f(b))} \&\& {\B(f(a),f(b)).}
	\arrow["{m_{\A(a,b)}}"', from=1-1, to=1-3]
	\arrow["{R\otimes f_{a,b}}"', from=1-1, to=3-1]
	\arrow["{m_{\B(f(a),f(b))}}", from=3-1, to=3-3]
	\arrow["{f_{a,b}}"', from=1-3, to=3-3]
\end{tikzcd}
\end{equation}
 Now we prove that this defines an $R$-$\mathrm{Mod}$-enriched functor.
 For graphical issues in the next diagrams, we will adopt the conventions that, for each $a,b\in\A$, we use the notation $\A_{a,b}$ instead of $\A(a,b)$, $\B_{fa,fb}$ instead of $\B(f(a),f(b))$, $m_{a,b}$ instead of $m_{\A(a,b)}$ and $m_{fa,fb}$ instead of $m_{\B(f(a),f(b))}$.
 
 To prove axiom \eqref{eqFirstDigramdefVenrichedFunctor} we need to show that the following diagram is commutative:
\begin{equation}
 \label{eqCorrezioneErrore2}   
\begin{tikzcd}[ampersand replacement=\&,sep=small]
	{R\otimes_{\Z}(\A_{a,b}\otimes_{R}\A_{b,c})} \&\& {R\otimes_{\Z}\A_{a,c}} \\
	{\A_{a,b}\otimes_{R}\A_{b,c}} \&\& {\A_{a,c}} \\
	\& {R\otimes_{\Z}(\B_{fa,fb}\otimes_{R}\B_{fb,fc})} \&\& {R\otimes_{\Z}\B_{fa,fc}} \\
	\& {\B_{fa,fb}\otimes_{R}\B_{fb,fc}} \&\& {\B_{fa,fc}.}
	\arrow[from=1-1, to=1-3]
	\arrow[from=1-1, to=2-1]
	\arrow[from=2-1, to=2-3]
	\arrow[from=1-3, to=2-3]
	\arrow[from=1-1, to=3-2]
	\arrow[from=2-1, to=4-2]
	\arrow[from=3-2, to=3-4]
	\arrow[from=1-3, to=3-4]
	\arrow[from=2-3, to=4-4]
	\arrow[from=4-2, to=4-4]
	\arrow[from=3-2, to=4-2]
	\arrow[from=3-4, to=4-4]
\end{tikzcd}
\end{equation}
We prove it by constructing the cube.
First, we build the bottom and right square.
Since the second component in $(f,\square)$ is a square in $Cat(\Ab)$, for each $a,b,c$ in $\A$ the following solid diagram is commutative:
\begin{equation}
\label{eqCorrezioneErrore3}    
\begin{tikzcd}[ampersand replacement=\&,sep=small]
	{\Z\otimes_{\Z}\A_{a,b}\otimes_{\Z}R\otimes_{\Z}\A_{b,c}} \&\& {\Z\otimes_{\Z}\B_{a,b}\otimes_{\Z}R\otimes_{\Z}\B_{b,c}} \\
	{R\otimes_{\Z}\A_{a,b}\otimes_{\Z}R\otimes_{\Z}\A_{b,c}} \&\& {R\otimes_{\Z}\B_{fa,fb}\otimes_{\Z}R\otimes_{\Z}\B_{fb,fc}} \\
	{\A_{a,b}\otimes_{\Z}\A_{b,c}} \&\& {\B_{fa,fb}\otimes_{\Z}\B_{fb,fc}} \\
	{\A_{a,b}\otimes_{R}\A_{b,c}} \&\& {\B_{fa,fb}\otimes_{R}\B_{fb,fc}} \\
	\& {R\otimes_{\Z}\A_{a,c}} \&\& {R\otimes_{\Z}\B_{fa,fc}} \\
	\& {\A_{a,c}} \&\& {\B_{fa,fc}.}
	\arrow["{u_{\A(a,b)}\otimes_{\Z}id}", shift left=2, from=1-1, to=2-1]
	\arrow["{id\otimes_{\Z}u_{\A(b,c)}}"', shift right=2, from=1-1, to=2-1]
	\arrow["{u_{\B(fa,fb)}\otimes_{\Z}id}", shift left=2, from=1-3, to=2-3]
	\arrow["{id\otimes_{\Z}u_{\B(fb,fc)}}"', shift right=2, from=1-3, to=2-3]
	\arrow["{m_{a,b}\otimes_{\Z} m_{b,c}}"', from=2-1, to=3-1]
	\arrow["{m_{fa,fb}\otimes_{\Z} m_{fb,fc}}"', from=2-3, to=3-3]
	\arrow["q"', from=3-1, to=4-1]
	\arrow["q"', from=3-3, to=4-3]
	\arrow["{f_{a,b}\otimes_{\Z}f_{b,c}}"{description}, from=3-1, to=3-3]
	\arrow["{R\otimes_{\Z}f_{a,b}\otimes_{\Z}R\otimes_{\Z}f_{b,c}}"', from=2-1, to=2-3]
	\arrow["\circ"{description}, from=3-1, to=6-2]
	\arrow["\circ"{description}, from=2-1, to=5-2]
	\arrow["{m_{a,c}}", from=5-2, to=6-2]
	\arrow[dashed, from=4-1, to=6-2]
	\arrow[dashed, from=4-1, to=4-3]
	\arrow[dashed, from=4-3, to=6-4]
	\arrow["\circ"{description}, from=3-3, to=6-4]
	\arrow["\circ"{description}, from=2-3, to=5-4]
	\arrow["{m_{fa,fc}}", from=5-4, to=6-4]
	\arrow["{f_{a,c}}"', from=6-2, to=6-4]
	\arrow["{R\otimes_{\Z}f_{a,c}}", from=5-2, to=5-4]
\end{tikzcd}
\end{equation}
In particular, the $q's$ are the coequalizers of the vertical diagrams above them: 
\[\begin{tikzcd}[ampersand replacement=\&]
	{\Z\otimes_{\Z}\A_{a,b}\otimes_{\Z}R\otimes_{\Z}\A_{b,c}} \& {R\otimes_{\Z}\A_{a,b}\otimes_{\Z}R\otimes_{\Z}\A_{b,c}} \& {\A_{a,b}\otimes_{\Z}\A_{b,c};}
	\arrow[shift left=2, from=1-1, to=1-2]
	\arrow[shift right=2, from=1-1, to=1-2]
	\arrow[from=1-2, to=1-3]
	\arrow["{m_{a,b}\otimes id }"',curve={height=18pt}, from=1-1, to=1-3]
	\arrow["{id\otimes m_{b,a}}",curve={height=-18pt}, from=1-1, to=1-3]
\end{tikzcd}\]
and 
\[\begin{tikzcd}[ampersand replacement=\&]
	{\Z\otimes_{\Z}\B_{fa,fb}\otimes_{\Z}R\otimes_{\Z}\B_{fb,fc}} \& {R\otimes_{\Z}\B_{fa,fb}\otimes_{\Z}R\otimes_{\Z}\B_{fb,fc}} \& {\B_{fa,fb}\otimes_{\Z}\B_{fb,fc};}
	\arrow[shift left=2, from=1-1, to=1-2]
	\arrow[shift right=2, from=1-1, to=1-2]
	\arrow[from=1-2, to=1-3]
	\arrow["{m_{fa,fb}\otimes id }"', curve={height=18pt}, from=1-1, to=1-3]
	\arrow["{id\otimes m_{fb,fa}}", curve={height=-18pt}, from=1-1, to=1-3]
\end{tikzcd}\]
see \eqref{eqCorrezioneErrore5}.
Using the universal property of coequalizers we obtain the dotted arrows and the diagram in \eqref{eqCorrezioneErrore3} is commutative by construction. This implies that both the lower than the right squares in \eqref{eqCorrezioneErrore2} are commutative.
The upper square of \eqref{eqCorrezioneErrore2} is commutative because it is the lower square tensorized via $R\otimes_{\Z}-$. The remaining of the diagram is computed as in \eqref{eqComposition} and so is canonically commutative.

To prove \eqref{eqSecondDigramdefVenrichedFunctor}, we need to prove that for each object $a\in\A$ the following diagram is commutative:
\[\begin{tikzcd}[ampersand replacement=\&,sep=small]
	\&\& {R\otimes_{\Z}(R\otimes_{\Z}\B(fa,fa))} \\
	{R\otimes_{\Z}R} \& {R\otimes_{\Z}(R\otimes_{\Z}\Z)} \& {R\otimes_{\Z}(R\otimes_{\Z}\A(a,a))} \& {R\otimes_{\Z}\A(a,a)} \& {R\otimes_{\Z}\B(fa,fa)} \\
	R \& {R\otimes_{\Z}\Z} \& {R\otimes_{\Z}\A(a,a)} \& {\A(a,a)} \& {\B(fa,fa).} \\
	\&\& {R\otimes_{\Z}\B(fa,fa)}
	\arrow[from=2-1, to=3-1]
	\arrow["\simeq"{description}, draw=none, from=2-1, to=2-2]
	\arrow["\simeq"{description}, draw=none, from=3-1, to=3-2]
	\arrow[from=2-2, to=3-2]
	\arrow[from=2-2, to=2-3]
	\arrow[from=3-2, to=3-3]
	\arrow[from=2-3, to=2-4]
	\arrow[from=3-3, to=3-4]
	\arrow[from=2-4, to=3-4]
	\arrow[from=2-4, to=2-5]
	\arrow[from=3-4, to=3-5]
	\arrow[from=2-5, to=3-5]
	\arrow[from=2-2, to=1-3]
	\arrow[from=3-2, to=4-3]
	\arrow[from=1-3, to=2-5]
	\arrow[from=4-3, to=3-5]
	\arrow[from=2-3, to=1-3]
	\arrow[from=3-3, to=4-3]
\end{tikzcd}\]
But this is commutative because it is composed of commutative diagrams. The central squares are commutative by definition of unit \eqref{eqUnit} and by the commutative of the ($a,a$)-component of the square \eqref{eqDajeChecisiamo}. The restarting diagrams are commutative by axioms \eqref{eqMorphismOfModuleassociativity} of \Cref{defCategoryOfLeftModuleObjects} (obviously, we take the ($a,a$)-component of the cube). 
We deduce that $\psi(f)$ is an $R\text{-}\mathrm{Mod}$-functor.

In turn, $\psi$ is a functor between categories because it preserves the composition.
Indeed, let $(f,\square_{f})$ and $(g,\square_{g})$ are two composable arrows in $\mathrm{LMod}_{\un{R}}(Cat(\Ab))$ represented by the commutative diagram:
\begin{equation}
\label{eqCorrezioneErrore4}
\begin{tikzcd}[ampersand replacement=\&]
	{\un{R}\otimes^{\Ab}\A} \& {\un{R}\otimes^{\Ab}\B} \& {\un{R}\otimes^{\Ab}\C} \\
	\A \& \B \& {\C.}
	\arrow["{\un{R}\otimes^{\Ab}f}"', from=1-1, to=1-2]
	\arrow["{\un{R}\otimes^{\Ab}g}"', from=1-2, to=1-3]
	\arrow[from=1-1, to=2-1]
	\arrow["f"', from=2-1, to=2-2]
	\arrow["g"', from=2-2, to=2-3]
	\arrow[from=1-3, to=2-3]
	\arrow["{\un{R}\otimes^{\Ab}gf}", curve={height=-18pt}, from=1-1, to=1-3]
	\arrow["gf"', curve={height=18pt}, from=2-1, to=2-3]
	\arrow[from=1-2, to=2-2]
\end{tikzcd}
\end{equation}
Let $a,b$ be two objects of $\A$.
The ($a,b$)-component of \eqref{eqCorrezioneErrore4} is the commutative diagram:  
\[\begin{tikzcd}[ampersand replacement=\&]
	{\un{R}\otimes^{\Ab}\A(a,b)} \& {\un{R}\otimes^{\Ab}\B(fa,fb)} \& {\un{R}\otimes^{\Ab}\C(gfa,gfb)} \\
	{\A(a,b)} \& {\B(fa,fb)} \& {\C(gfa,gfb);}
	\arrow["{\un{R}\otimes^{\Ab}f_{a,b}}"', from=1-1, to=1-2]
	\arrow["{\un{R}\otimes^{\Ab}g_{fa,fb}}"', from=1-2, to=1-3]
	\arrow[from=1-1, to=2-1]
	\arrow["{f_{a,b}}"', from=2-1, to=2-2]
	\arrow["{g_{fa,fb}}"', from=2-2, to=2-3]
	\arrow[from=1-3, to=2-3]
	\arrow["{\un{R}\otimes^{\Ab}gf_{a,b}}", curve={height=-18pt}, from=1-1, to=1-3]
	\arrow["{gf_{a,b}}"', curve={height=18pt}, from=2-1, to=2-3]
	\arrow[from=1-2, to=2-2]
\end{tikzcd}\]
which implies that $\psi$ respect the composition.

Now we prove that $\psi$ is an equivalence of categories.
By construction (it is a kind of identity on the objects) $\psi$ is essentially subjective.
We just have to prove that it is fully faithful. 
We will define the inverse of $\psi_{\A,\B}$, we will denote it by $\gamma$.
Let $h:\A\to\B$ be an arrow in $Cat(R\text{-}\mathrm{Mod})$ and let $U: Cat(R\text{-}\mathrm{Mod})\to Cat(\Ab)$ the canonically forgetful (i.e. the change of enrichment along the lax monoidal functor restriction of scalar $u_{R}^*:R\text{-}\mathrm{Mod}\to \Ab$). 
We define a $\un{R}$-equivariant arrow $\gamma(h):=(\hat{h},\square)$. 
As the first component we chose the arrow in $Cat(\Ab)$ \[\hat{h}:=U(h):U(\A)\to U(\B);\] 
the second component is the square: 
\[\begin{tikzcd}[ampersand replacement=\&]
	{\un{R}\otimes_{\Z}\A} \&\& U\A \\
	\\
	{\un{R}\otimes_{\Z}\B} \&\& {U\B;}
	\arrow[from=1-1, to=1-3]
	\arrow["{R\otimes_{\Z}Uh}"', from=1-1, to=3-1]
	\arrow[from=3-1, to=3-3]
	\arrow["{Uh=\hat{h}}", from=1-3, to=3-3]
\end{tikzcd}\]
whose ($a,b$)-component is \eqref{eqCorrezioneErrore6}. 
The proof that square is well-define (i.e. it is composed by $Cat(\Ab)$-enriched functor) is exactly the commutativity of the diagram \eqref{eqCorrezioneErrore3} (including the dotted arrows). 
Practically, we used another time the universal properties of the coequalizer but in the opposite direction.

By construction, for each $f\in \mathrm{LMod}_{\un{R}}(Cat(\Ab))$ and $h\in Cat(R\text{-}\mathrm{Mod})$ the equations $\psi_{\A,\B}(\gamma({h}))=h$ and $\gamma(\psi_{\A,\B}(f))=f$ hold, this implies that $\psi$ is fully faithful.
So we have proved the theorem.
\end{proof}

We need the following lemma to prove the second equivalence in \Cref{sumtheorem}.

\begin{lemma}
\label{lmAdInduces2ad}
Let $(L\dashv R:\A\to \B, \eta: \text{id} \Rightarrow RL, \epsilon: LR \Rightarrow \text{id})$ be an adjunction, then there exists an adjunction between the arrows categories
$(L^{\mathbb{2}}\dashv R^{\mathbb{2}}:\A^{\mathbb{2}}\to \B^{\mathbb{2}}, \eta^{\mathbb{2}}: \text{id} \Rightarrow RL, \epsilon^{\mathbb{2}}: L^{\mathbb{2}}R^{\mathbb{2}} \Rightarrow \text{id}^{\mathbb{2}})$
where for each arrow $f:a\to a'$ of $\A$ and $g:b\to b'$ of $\B$ the unit and the counit components are the following squares:

\[
\begin{tikzcd}[ampersand replacement=\&]
	{\eta^{\mathbb{2}}_f:=} \& a \&\& RLa \& {\epsilon_{g}^{\mathbb{2}}:=} \& LRb \&\& b \\
	\& {a'} \&\& {RLa';} \&\& {LRb'} \&\& {b'.}
	\arrow["{\eta_a}"', from=1-2, to=1-4]
	\arrow["f", from=1-2, to=2-2]
	\arrow["{\eta_{a'}}", from=2-2, to=2-4]
	\arrow["{R^\mathbb{2}L^\mathbb{2}f}"', from=1-4, to=2-4]
	\arrow["{L^\mathbb{2}R^\mathbb{2}g}"', from=1-6, to=2-6]
	\arrow["{\eta_{b'}}", from=2-6, to=2-8]
	\arrow["{\eta_b}"', from=1-6, to=1-8]
	\arrow["g", from=1-8, to=2-8]
\end{tikzcd}
\]
    
\end{lemma}

\begin{proof}
We denote by $2Cat$ the $Cat$-enriched category (or, equivalently, a $2$-category) of categories.

The internal-hom is a $Cat$-enriched functor (or a $2$-functor) 
\[(-)^\mathbb{2}:2Cat\to 2Cat: \A\to \A^{\mathbb{2}}.\]

Since any $Cat$-enriched functor preserves (internal) adjunction, \cite[Exercise B.3.iv]{RiehlElements}, it sends the adjunction $L\dashv R$ into the adjunction $L^{2}\dashv R^{2}$. So, the thesis.
    
\end{proof}

So we can prove the second main result.

\begin{theorem}
\label{thcategorialversionAction}
Let $R$ be a commutative ring with unit. There is an equivalence of categories
\begin{equation*}
   Fun^{Cat(\Ab)}(\2un{R},Cat(\Ab))\simeq \mathrm{LMod}_{\underline{R}}(Cat(\Ab)).
\end{equation*}
\end{theorem}

\begin{proof}
Let us start by defining a functor
\[\gamma: Fun^{Cat(\Ab)}(\2un{R},Cat(\Ab))\to \mathrm{LMod}_{\underline{R}}(Cat(\Ab)). \]

An object $F$ of $Fun^{Cat(\Ab)}(\2un{R},Cat(\Ab))$ is defined by a morphism in $Set$ 
\begin{equation}
\label{eqRealFirstComponent}
    obF:\{*\}\simeq ob(\2un{R}) \to ob(Cat(\Ab)): *\mapsto F(*); 
\end{equation}    
and by a functor of monoids in $Cat(\Ab)$
\begin{equation}
\label{eqRealSecondComponent}    
F_{**}:\2un{R}(*,*)\simeq \underline{R}\to F(*)^{F(*)},
\end{equation}
where $F(*)^{F(*)}$ is a commutative monoid with the composition as multiplication and the identity as unit.

In turn, $F_{**}$ is defined by a morphism in $Set$
\begin{equation}
\label{eqFirstComponent}
    ob{F_{**}}: \{*\}\simeq ob(\underline{R})\to ob(F(*)^{F(*)}):*\mapsto F_{**}(*); 
\end{equation}
and a morphism of rings, equivalently a morphism of monoids in $\Ab$  
\begin{equation}
\label{eqSecondComponent}    
F_{**,**}:\underline{R}(*,*)\simeq R\to F(*)^{F(*)}(F_{**}(*),F_{**}(*)).
\end{equation}

Now we define $\gamma$. $\gamma$ sends an object $F:\2un{R}\to Cat(\Ab)$ of $Fun^{Cat(\Ab)}(\2un{R},Cat(\Ab))$ into the object $(F(*), m_{F(*)}:\underline{R}\otimes^{\Ab} F(*)\to F(*))$ of $\mathrm{LMod}_{\underline{R}}(Cat(\Ab))$, whose first component is defined in (\ref{eqRealFirstComponent}) and whose second component is the $\Ab$-enriched functor 
\[m_{F(*)}:\underline{R}\otimes^{\Ab} F(*)\to F(*);\]
 which is the adjoint-correspondent of the $\Ab$-enriched functor \eqref{eqRealSecondComponent} via the tensor-internal-hom adjunction $-\otimes F(*)\dashv (-)^{F(*)}:Cat(\Ab)\to Cat(\Ab)$.

Furthermore, for each pair $F,G$ of objects of $ Fun^{Cat(\Ab)}(\2un{R},Cat(\Ab))$, we define the morphism of $Set$ \[\gamma_{F,G}:Fun^{Cat(\Ab)}(\2un{R},Cat(\Ab))(F,G)\to \mathrm{LMod}_{\underline{R}}(Cat(\Ab))(\gamma(F),\gamma(G)).\]

Let $\alpha:F\Rightarrow G$ be an object of $Fun^{Cat(\Ab)}(\2un{R},Cat(\Ab))(F,G)$. Since $\2un{R}$ has only one object, $\alpha$ consists of only one morphism $\alpha:\underline{\Z}\to Cat(\Ab)(F(*),G(*)) $ in $Cat(\Ab)$ such that only one square is commutative: 
\begin{equation}
\label{eqEnrichednaturalTransformation}  
\begin{tikzcd}[ampersand replacement=\&]
	{\2un{R}(*,*)\simeq\un{R}} \&\& {F(*)^{F(*)}} \\
	\\
	{G(*)^{G(*)}} \&\& {G(*)^{F(*)};}
	\arrow["{F_{**}}", from=1-1, to=1-3]
	\arrow["{G_{**}}"', from=1-1, to=3-1]
	\arrow["{\alpha^{*}}"', from=3-1, to=3-3]
	\arrow["{\alpha_{*}}", from=1-3, to=3-3]
\end{tikzcd}
\end{equation}
see \Cref{defEnrichedNaturalTransformation}. Above we made an abuse: we denoted by $\alpha:F(*)\to G(*)$ instead that $\alpha_{*}$ the only one component of $\alpha:F\Rightarrow G$. We chose to make this abuse so as not to confuse the component of $\alpha$ and postcomposition with $\alpha$ (\Cref{notComposition}).

We have the two tensor-internal-hom adjunction \[(-\otimes F(*)\dashv (-)^{F(*)}: Cat(\Ab)\to Cat(\Ab), \epsilon^{F},\eta^{F});\] and  \[(-\otimes G(*)\dashv (-)^{G(*)}: Cat(\Ab)\to Cat(\Ab), \epsilon^{G},\eta^{G});\]  in both the triple the last two components are the respective units and counits.

Every adjunction induces an adjunction between the categories of arrows, see \Cref{lmAdInduces2ad}, the adjunction $-\otimes F(*)\dashv (-)^{F(*)}: Cat(\Ab)\to Cat(\Ab)$ induces the adjunction \[((-\otimes F(*))^{\mathbb{2}}\dashv ((-)^{F(*)})^{\mathbb{2}}:(Cat(\Ab))^{\mathbb{2}}\to Cat(\Ab), \epsilon^{F,\mathbb{2}},\eta^{F,\mathbb{2}}).\]
Using the universal property of this last adjunction on (\ref{eqEnrichednaturalTransformation}) we obtain the following commutative diagram:
\begin{equation}
    \label{eqFirtsTerribleDiagram}
    \begin{tikzcd}[ampersand replacement=\&]
	{\un{R}\otimes^{\Ab}F(*)} \&\&\&\& {F(*)^{F(*)}\otimes^{\Ab}F(*)} \&\& {F(*)} \\
	{\un{R}\otimes^{\Ab}F(*)} \&\& {G(*)^{G(*)}\otimes^{\Ab}F(*)} \&\& {G(*)^{F(*)}\otimes^{\Ab}F(*)} \&\& {G(*).}
	\arrow[shift left=2, no head, from=1-1, to=2-1]
	\arrow["{ G_{**}\otimes^{\Ab}F(*)}"', from=2-1, to=2-3]
	\arrow["{\alpha^{*}\otimes^{\Ab}F(*)}"', from=2-3, to=2-5]
	\arrow["{ F_{**}\otimes^{\Ab}F(*)}", from=1-1, to=1-5]
	\arrow["{\alpha_{*}\otimes^{\Ab}F(*)}", from=1-5, to=2-5]
	\arrow["{\epsilon_{G(*)}^{F}}", from=2-5, to=2-7]
	\arrow["{\epsilon_{F(*)}^{F}}", from=1-5, to=1-7]
	\arrow["\alpha"', from=1-7, to=2-7]
	\arrow[shift right=2, no head, from=1-1, to=2-1]
\end{tikzcd}
\end{equation}
Now, we prove that the bottom horizontal line factor through \[\underline{R}\otimes F(*)\xrightarrow{\underline{R}\otimes \alpha} \underline{R}\otimes G(*).\]
 We archive it to prove that the following diagram is commutative:
\begin{equation}
\label{eqSecondTerribleDiagram} 
\begin{tikzcd}[ampersand replacement=\&]
	{\un{R}\otimes^{\Ab}F(*)} \&\& {G(*)^{G(*)}\otimes^{\Ab}F(*)} \&\& {G(*)^{F(*)}\otimes^{\Ab}F(*)} \&\& {G(*)} \\
	\\
	{\un{R}\otimes^{\Ab}G(*)} \&\& {G(*)^{G(*)}\otimes^{\Ab}G(*)} \&\&\&\& {G(*).}
	\arrow["{ G_{**}\otimes^{\Ab}F(*)}"', from=1-1, to=1-3]
	\arrow["{\alpha_{*}\otimes F(*)}"', from=1-3, to=1-5]
	\arrow["{\epsilon_{G(*)}^{F}}", from=1-5, to=1-7]
	\arrow["{\un{R}\otimes \alpha}", from=1-1, to=3-1]
	\arrow["{G_{**}\otimes^{\Ab}G(*)}", from=3-1, to=3-3]
	\arrow["{id\otimes^{\Ab}\alpha}"', from=1-3, to=3-3]
	\arrow["{\epsilon_{G(*)}^{G}}", from=3-3, to=3-7]
	\arrow[shift right, no head, from=1-7, to=3-7]
	\arrow[shift left, no head, from=1-7, to=3-7]
\end{tikzcd}
\end{equation}

Indeed, the right rectangle commutes for the dinaturality of the counit, see \cite[\S 1.1.1]{loregian2021co} and the left rectangle commutes because this square is the image of the Leibniz tensor, \Cref{defLeibnizConstruction}, 
\[\hat{-\otimes-}:Cat(\Ab)^\mathbb{2}\times Cat(\Ab)^\mathbb{2}\to Cat(\Ab)^\mathbb{2};\]
of the pair: 
\[\begin{tikzcd}[ampersand replacement=\&]
	{\un{R}} \& {\un{Z}} \& {\un{Z}} \\
	{G(*)^{G(*)}} \& {F(*)} \& {G(*);}
	\arrow[from=1-1, to=2-1]
	\arrow[shift left, no head, from=1-2, to=1-3]
	\arrow["\emptyset"', from=1-2, to=2-2]
	\arrow["\alpha"', from=2-2, to=2-3]
	\arrow["\emptyset", from=1-3, to=2-3]
	\arrow[shift right, no head, from=1-2, to=1-3]
\end{tikzcd}\]
where $\emptyset$ denotes the unique arrows from the initial object of $Cat(\Ab)$.

So, we define $\gamma_{F,G}(\alpha)\in \mathrm{LMod}_{\underline{R}}(Cat(\Ab))(\gamma(F),\gamma(G)) $ as the commutative square obtained by pasting together (\ref{eqFirtsTerribleDiagram}) and \eqref{eqSecondTerribleDiagram} along the two equal horizontal lines:   
\[\begin{tikzcd}[ampersand replacement=\&]
	{\un{R}\otimes^{\Ab}F(*)} \&\& {F(*)} \\
	{\un{R}\otimes^{\Ab}G(*)} \&\& {G(*).}
	\arrow["{\un{R}\otimes^{\Ab}\alpha}"', from=1-1, to=2-1]
	\arrow[from=2-1, to=2-3]
	\arrow[from=1-1, to=1-3]
	\arrow["\alpha", from=1-3, to=2-3]
\end{tikzcd}\]
Now, we prove that $\gamma$ satisfies the functoriality axioms in \Cref{defEnrichedFunctor}.
It is trivial that $\gamma$ sends identity morphism in identity morphism.
We only have to prove that $\gamma$ respects the composition.
Let $\alpha:F\Rightarrow G$ and $\beta: G\Rightarrow H$ be two arrows of $Fun^{Cat(\Ab)}(\2un{R},Cat(\Ab))$, we must to prove that $\gamma (\alpha)\circ \gamma (\beta)= \gamma(\beta\circ \alpha)$. This is equivalent to proving that the perimeters of the following commutative diagram (which is $\gamma (\alpha)\circ \gamma (\beta)$): 
\[\begin{tikzcd}[ampersand replacement=\&]
	{\un{R}\otimes^{\Ab}F(*)} \&\&\&\& {F(*)\otimes^{\Ab}F(*)^{F(*)}} \&\& {F(*)} \\
	{\un{R}\otimes^{\Ab}F(*)} \&\& \bullet \&\& \bullet \&\& {G(*)} \\
	{\un{R}\otimes^{\Ab}G(*)} \&\& \bullet \&\&\&\& {G(*)} \\
	{\un{R}\otimes^{\Ab}G(*)} \&\& \bullet \&\& \bullet \&\& {H(*)} \\
	{\un{R}\otimes^{\Ab}H(*)} \&\&\&\& {H(*)\otimes^{\Ab}H(*)^{H(*)}} \&\& {H(*);}
	\arrow[shift left, no head, from=1-1, to=2-1]
	\arrow[from=2-1, to=2-3]
	\arrow[from=2-3, to=2-5]
	\arrow[from=2-5, to=2-7]
	\arrow["\alpha"', from=1-7, to=2-7]
	\arrow["{\un{R}\otimes^{\Ab}\alpha}"', from=2-1, to=3-1]
	\arrow[from=3-1, to=3-3]
	\arrow[from=2-3, to=3-3]
	\arrow[from=3-3, to=3-7]
	\arrow[shift right, no head, from=2-7, to=3-7]
	\arrow[shift left, no head, from=3-1, to=4-1]
	\arrow["{\un{R}\otimes^{\Ab}\beta}"', from=4-1, to=5-1]
	\arrow["{F(*)\otimes^{\Ab}F_{**}}"', from=1-1, to=1-5]
	\arrow["{\epsilon^{F}_{F(*)}}"', from=1-5, to=1-7]
	\arrow[from=3-3, to=4-5]
	\arrow[from=1-5, to=2-5]
	\arrow[from=4-1, to=4-3]
	\arrow[from=4-3, to=4-5]
	\arrow[from=4-5, to=4-7]
	\arrow["\beta"', from=3-7, to=4-7]
	\arrow["{H(*)\otimes^{\Ab}H_{**}}"', from=5-1, to=5-5]
	\arrow["{\epsilon_{H(*)}^{H}}"', from=5-5, to=5-7]
	\arrow[shift right, no head, from=4-7, to=5-7]
	\arrow[shift right, no head, from=1-1, to=2-1]
	\arrow[shift right, no head, from=3-1, to=4-1]
	\arrow[shift left, no head, from=2-7, to=3-7]
	\arrow[shift left, no head, from=4-7, to=5-7]
\end{tikzcd}\]

and the perimeters of the following commutative diagram (which is $\gamma(\beta\circ \alpha)$):

\[\begin{tikzcd}[ampersand replacement=\&]
	{\un{R}\otimes^{\Ab}F(*)} \&\&\&\& {F(*)\otimes^{\Ab}F(*)^{F(*)}} \&\& {F(*)} \\
	{\un{R}\otimes^{\Ab}F(*)} \&\& \bullet \&\& \bullet \&\& {G(*)} \\
	{\un{R}\otimes^{\Ab}H(*)} \&\& {H(*)^{H(*)}\otimes^{\Ab}H(*)} \&\&\&\& {G(*);}
	\arrow[shift left, no head, from=1-1, to=2-1]
	\arrow[from=2-1, to=2-3]
	\arrow[from=2-3, to=2-5]
	\arrow[from=2-5, to=2-7]
	\arrow["\beta\circ\alpha"', from=1-7, to=2-7]
	\arrow["{\un{R}\otimes^{\Ab}\beta\circ\alpha}", from=2-1, to=3-1]
	\arrow["{H(*)\otimes^{\Ab}H_{**}}"', from=3-1, to=3-3]
	\arrow[from=2-3, to=3-3]
	\arrow["{\epsilon_{H(*)}^{H}}"', from=3-3, to=3-7]
	\arrow[shift right, no head, from=2-7, to=3-7]
	\arrow["{F(*)\otimes^{\Ab}F_{**}}"', from=1-1, to=1-5]
	\arrow["{\epsilon^{F}_{F(*)}}"', from=1-5, to=1-7]
	\arrow[from=1-5, to=2-5]
	\arrow[shift right, no head, from=1-1, to=2-1]
	\arrow[shift left, no head, from=2-7, to=3-7]
\end{tikzcd}\]
are equal. For simplicity, we have replaced information no longer necessary in the diagrams above with a dot. This holds because the only sides that can be different are the left vertical ones, but the two left vertical sides are equal for the funtoriality of the bifunctor $-\otimes^{\Ab}-$.
We conclude that $\gamma$ is a functor.

Now we prove that it is an equivalence of categories.
$\gamma$ is essentially surjective. Indeed, let $(\A,m_{\A}:\underline{R}\otimes \A\to \A)$ be an object of $\mathrm{LMod}_{\underline{R}}(Cat(\Ab))$, we define the $Cat(\Ab)$-enriched 
functor $F:\2un{R}\to Cat^{Cat(\Ab)}$ on the set of objects as:
\[ob(F):\{*\}=ob(\2un{R})\to ob(Cat(\Ab)):*\mapsto \A;\]
and on the arrows by the arrow of rings (or monoids)   
\[F_{**}:\2un{R}(*,*)\simeq\underline{R}\to \A^{\A};\]
which is the adjoint-correspondent of $m_{\A}$ via the adjunction $-\otimes^{\Ab}\A\dashv (-)^{\A}$. 
By construction of $\gamma$, we obtain that $\gamma(F)\simeq (\A,m_{\A}:\underline{R}\otimes \A\to \A)$.

It is also fully faithful.
We will define the inverse of $\gamma_{F,G}$, we will denote it by $\psi$.

Let $f$ be an arrow in $\mathrm{LMod}_{\underline{R}}(Cat(\Ab))(\gamma(F),\gamma(G))$ which is represented by the following square:

\begin{equation}
\label{eqCenaRitardo}    
\begin{tikzcd}[ampersand replacement=\&]
	{\un{R}\otimes^{\Ab} F(*)} \&\& {F(*)} \\
	{\un{R}\otimes^{\Ab} G(*)} \&\& {G(*).}
	\arrow["{m_{F(*)}}"', from=1-1, to=1-3]
	\arrow["{\un{R}\otimes^{\Ab} f}", from=1-1, to=2-1]
	\arrow["{m_{G(*)}}", from=2-1, to=2-3]
	\arrow["f"', from=1-3, to=2-3]
\end{tikzcd}
\end{equation}

We define $\psi(f)$ as the following commutative square:
\[\begin{tikzcd}[ampersand replacement=\&]
	{\un{R}} \&\& {(\un{R}\otimes^{\Ab} F(*))^{F(*)}} \&\& {F(*)^{F(*)}} \\
	\&\& {(\un{R}\otimes^{\Ab} G(*))^{F(*)}} \&\& {G(*)^{F(*)}} \\
	{\un{R}} \&\& {(\un{R}\otimes^{\Ab} G(*))^{G(*)}} \&\& {G(*)^{G(*)}.}
	\arrow["{m_{F(*)}^{F(*)}}"', from=1-3, to=1-5]
	\arrow["{\un{R}\otimes^{\Ab} f}", from=1-3, to=2-3]
	\arrow["{m_{G(*)}^{F(*)}}", from=2-3, to=2-5]
	\arrow["{f^{F(*)}}"', from=1-5, to=2-5]
	\arrow["{\eta_{\un{R}}^{F(*)}}"', from=1-1, to=1-3]
	\arrow[shift left, no head, from=1-1, to=3-1]
	\arrow["{\eta_{\un{R}}^{G(*)}}", from=3-1, to=3-3]
	\arrow["{(\un{R}\otimes^{\Ab} G(*))^{f}}", from=3-3, to=2-3]
	\arrow["{m_{G(*)}^{G(*)}}", from=3-3, to=3-5]
	\arrow["{G(*)^{f}}", from=3-5, to=2-5]
	\arrow[shift right, no head, from=1-1, to=3-1]
\end{tikzcd}\]
The above diagram is commutative because is composed of commutative diagrams. The left rectangle is commutative by the dinaturality of the unit $\eta$, the right-top rectangle is commutative because it is the image of \eqref{eqCenaRitardo} via the functor $(-)^{F(*)}$, and the right-bottom rectangle is commutative because it is the image of the pair:
\[\begin{tikzcd}[ampersand replacement=\&]
	{\un{R}\otimes^{\Ab} G(*)} \& {G(*)} \& {F(*)} \\
	{\un{R}\otimes^{\Ab} G(*)} \& {G(*)} \& {G(*);}
	\arrow[shift left, no head, from=1-1, to=2-1]
	\arrow["{m_{G(*)}}"', from=1-1, to=1-2]
	\arrow[shift left, no head, from=1-2, to=2-2]
	\arrow["{m_{G(*)}}", from=2-1, to=2-2]
	\arrow[shift right, no head, from=1-1, to=2-1]
	\arrow[shift right, no head, from=1-2, to=2-2]
	\arrow[from=1-3, to=2-3]
\end{tikzcd}\]
via the Leibniz internal-hom, \Cref{defLeibnizConstruction},
\[\hat{(-)^{(-)}}:(Cat(\Ab))^{\mathbb{2}}\times( Cat(\Ab))^{\mathbb{2}}\to (Cat(\Ab))^{\mathbb{2}}.  \]

Finally, we prove that $\psi$ is the inverse of $\gamma_{F,G}$ showing that the following equations holds:
\begin{equation}
\label{eqPenultimissimissima}
    \psi\circ \gamma_{F,G}= id;
\end{equation}
and
\begin{equation}
\label{eqUltimissimissimissima}    
\gamma_{F,G}\circ \psi= id.  
\end{equation}

Let $f$ be an arrow in $\mathrm{LMod}_{\un{R}}(Cat(Ab))(\gamma (F ), \gamma (G))$ represented by the square \eqref{eqCenaRitardo}.
The arrow $\gamma_{F,G}\circ \psi(f)$ is represented by the following square:
\begin{equation}
\label{eqFinalmenteQualcosaQuasiFinito}
\begin{tikzcd}[ampersand replacement=\&,sep=small]
	{\un{R}\otimes^{\Ab}F(*)} \&\& {(\un{R}\otimes^{\Ab}F(*)^{F(*)})\otimes^{\Ab}F(*)} \&\& {F(*)^{F(*)}\otimes^{\Ab}F(*)} \&\& {F(*)} \\
	{\un{R}\otimes^{\Ab}F(*)} \&\& \bullet \&\& \bullet \& \bullet \& {G(*)} \\
	{\un{R}\otimes^{\Ab}G(*)} \&\& {(\un{R}\otimes^{\Ab}G(*)^{G(*)})\otimes^{\Ab}G(*)} \&\& {G(*)^{G(*)}\otimes^{\Ab}G(*)} \&\& {G(*).}
	\arrow[shift right, no head, from=1-1, to=2-1]
	\arrow[shift left, no head, from=1-1, to=2-1]
	\arrow["{\un{R}\otimes^{\Ab}f}"', from=2-1, to=3-1]
	\arrow["{\;\eta^{F}_{\un{R}}\otimes^{\Ab}F(*)}", from=1-1, to=1-3]
	\arrow["{m^{F(*)}_{F(*)}\otimes F(*)}", from=1-3, to=1-5]
	\arrow["{\epsilon^{F}_{F(*)}}", from=1-5, to=1-7]
	\arrow["f"', from=1-7, to=2-7]
	\arrow[from=2-6, to=2-7]
	\arrow[from=1-5, to=2-6]
	\arrow[shift left, no head, from=2-7, to=3-7]
	\arrow[shift right, no head, from=2-7, to=3-7]
	\arrow["{\epsilon^{G}_{G(*)}}"', from=3-5, to=3-7]
	\arrow[from=2-5, to=2-6]
	\arrow["{\;\eta^{G}_{\un{R}}\otimes^{\Ab}G(*)}"', from=3-1, to=3-3]
	\arrow["\;{m^{G(*)}_{G(*)}\otimes^{\Ab} G(*)}", from=3-3, to=3-5]
	\arrow[from=2-5, to=3-5]
	\arrow[from=2-1, to=2-3]
	\arrow[from=2-3, to=2-5]
\end{tikzcd}
\end{equation}

For simplicity, we have replaced information no longer necessary in the diagrams above with a dot.
If we prove that the bottom and the upper horizontal lines in \eqref{eqFinalmenteQualcosaQuasiFinito} are the same of $m_{F(*)}$ and $m_{G(*)}$ respectively, we have proved \eqref{eqUltimissimissimissima}.
The upper line is the arrow $m_{F(*)}:\un{R}\otimes^{\Ab}F(*)\to F(*)$ after applying the twice the adjoint-correspondence
\begin{equation}
\label{eqDAJEEEEEEEEEE}
\Hom_{Cat(\Ab)}(-\otimes^{\Ab}F(*),-)\simeq \Hom_{Cat(\Ab)}(-,(-)^{F(*)}).
\end{equation}
Indeed, after one application we obtain the arrow
\[R\xrightarrow{\eta_{\un{R}}^{F(*)}}\un{R}\otimes^{\Ab}F(*)^{F(*)}\xrightarrow{m_{F(*)}^{F(*)}}F(*)^{F(*)};\]
and, after the second application, we obtain the upper arrow in \eqref{eqFinalmenteQualcosaQuasiFinito}. 
Since \eqref{eqDAJEEEEEEEEEE} is a bijection, we obtain that the upper arrow is the same as $m_{F(*)}$.
The case with $m_{G(*)}$ is similar but with the adjoint-correspondence 
\begin{equation}
\label{eqVAdaVialCu}
    \Hom_{Cat(\Ab)}(-\otimes^{\Ab}G(*),-)\simeq \Hom_{Cat(\Ab)}(-,(-)^{G(*)}).
\end{equation}
So, \eqref{eqUltimissimissimissima} holds.

Let $\alpha:F\Rightarrow G$ be an object of $Fun^{Cat(\Ab)}(\2un{R},Cat(\Ab))(F,G)$ as in \eqref{eqEnrichednaturalTransformation}.
The natural transformation $\gamma \psi_{F,G}(\alpha)$ is represented by the following commutative diagram:
\begin{equation}
\label{eqUltimaDEllArticolo}
\begin{tikzcd}[ampersand replacement=\&]
	{\un{R}} \&\& {(\un{R}\otimes^{\Ab}F(*))^{F(*)}} \&\& {(F(*)^{F(*)}\otimes^{\Ab}F(*))^{F(*)}} \&\& {F(*)^{F(*)}} \\
	\&\& \bullet \&\&\&\& {G(*)^{F(*)}} \\
	{\un{R}} \&\& {(\un{R}\otimes^{\Ab}G(*))^{G(*)}} \&\& {(F(*)^{F(*)}\otimes^{\Ab}F(*))^{F(*)}} \&\& {F(*)^{F(*)}.}
	\arrow["{\eta_{\un{R}}^{F}}", from=1-1, to=1-3]
	\arrow["{\; \;( F_{**}\otimes^{\Ab}F(*))^{F(*)}}", from=1-3, to=1-5]
	\arrow["{(\epsilon_{F(*)}^{F})^{F(*)}}", from=1-5, to=1-7]
	\arrow["{\alpha^{F(*)}}"', from=1-7, to=2-7]
	\arrow["{(\epsilon_{G(*)}^{G})^{G(*)}}", from=3-5, to=3-7]
	\arrow["{\;\; \;( G_{**}\otimes^{\Ab}G(*))^{G(*)}}", from=3-3, to=3-5]
	\arrow["{\eta_{\un{R}}^{G}}", from=3-1, to=3-3]
	\arrow[shift left, no head, from=1-1, to=3-1]
	\arrow[shift right, no head, from=1-1, to=3-1]
	\arrow["{G(*)^{\alpha}}", from=3-7, to=2-7]
	\arrow[from=1-3, to=2-3]
	\arrow[from=3-3, to=2-3]
	\arrow[from=2-3, to=2-7]
\end{tikzcd}
\end{equation}
For simplicity, we have replaced information no longer necessary in the diagrams above with a dot.
If we establish that the upper and bottom horizontal lines correspond to $F_{**}$ and $G_{**}$ respectively, we thereby prove \eqref{eqPenultimissimissima}.
The upper line is the arrow $F_{**}$ after applying twice the adjoint-correspondence \eqref{eqDAJEEEEEEEEEE}.
Since \eqref{eqDAJEEEEEEEEEE} is a bijection, we obtain that the upper arrow is the same as $m_{F(*)}$.
The case with $G_{**}$ is similar but with the adjoint-correspondence \eqref{eqVAdaVialCu}.

So, we proved the theorem.
\end{proof}


\begin{remark}
\label{rmkGeneralizzation}
There are some obvious generalizations of \Cref{thcategorialversionEnrichedModareMod}. One is to use a general symmetric monoidal category $\mathcal{V}$ instead of $\Ab$. Since the monoidal structure in $\mathcal{V}$ is symmetric, the category $CMon(\mathcal{V})$ has a monoidal structure such that the forgetful functor $CMon(\mathcal{V}) \to \mathcal{V}$ is monoidal. In particular, the underline functor is monoidal in this case.
Probably the results of this paper do not hold for each symmetric monoidal category $\V$. For example in our proof of \Cref{thcategorialversionAction} is crucial that the Leibniz tensor and internal-hom are well-defined. Thus, we expect that $\V$ must have enough properties to enable us to define the Leibniz construction in $Cat(\V)$.  
We think that the above generalization may work, but we are not sure. For example, there may be some size drawbacks. Adding some control on the colimits in $\mathcal{V}$ it is probably possible to overcome the size problem, e.g., if $\mathcal{V}$ is a locally presentable category and $\otimes_{\mathcal{V}}$ preserves colimits.
However, this generalization does not concern this paper.
\end{remark}


In the last part of this article, we will find a specific case of the \Cref{thcategorialversionEnrichedModareMod}. 
We spend some words on \eqref{eqAndreaBovo} and we will find the real \enquote{one object version} of \Cref{thcategorialversionEnrichedModareMod}, see \Cref{thRealOneObjectVersion}.   

Let $Ring$ be the ordinary category of rings with a unit. This category has a symmetric monoidal structure $(Ring, \otimes_{\mathbb{Z}}, \mathbb{Z})$ whose second component is the tensor product of rings (which is the same as the tensor of commutative groups). 

The underline functor $\underline{(-)}:Ring\simeq Mon(\Ab)\to Cat(\Ab)$ is fully faithful and monoidal. Indeed, by definition, the equivalence of $\Ab$-enriched categories holds:
\[\underline{M} \otimes^{\Ab} \underline{M} \simeq \underline{M\otimes_{\mathbb{Z}}M}. \]

Then we define $\underline{Ring} \subseteq Cat(\Ab)$ as the full subcategory generated by the image of the underline functor $\underline{(-)}: Ring \to Cat(\Ab): M \mapsto \underline{M} $ and, since $\underline{(-)}$ is a monoidal functor, it inherits the structure of a symmetric monoidal category.

Finally, we can state the \enquote{with one object version} of \Cref{thcategorialversionEnrichedModareMod}.

\begin{theorem}[One object version]
\label{thRealOneObjectVersion}
Let $R$ be a commutative ring with a unit. There exists a chain of equivalences of categories:
\begin{equation}
\begin{split}
     \mathrm{LMod}_{R}(Ring) & \simeq \mathrm{LMod}_{\underline{R}}(\underline{Ring}) \\
     & \simeq Fun^{Cat(\Ab)}(\underline{R},\underline{Ring}) \\
     & \simeq Fun^{\underline{Ring}}(\underline{R},\underline{Ring}).
\end{split}
\end{equation}
\end{theorem}

\begin{proof}
The same as \Cref{thcategorialversionAction} for the second equivalence. The first follows from the fact that the underline functor $\underline{(-)}: Ring\simeq Mon(\Ab)\to Cat(\Ab)$ is fully faithful and monoidal. The last is canonical.
\end{proof}

\vskip 1cm

 \printbibliography 
\end{document}